\documentclass[fleqn,a4paper,10pt]{amsart}

\usepackage{graphicx}
\usepackage{tikz}
\usetikzlibrary{shapes,decorations.pathmorphing}
\usepackage{amsmath}
\usepackage{amssymb,amsthm}
\usepackage{comment}
\usepackage{hyperref}
\usepackage{enumerate}

\theoremstyle{definition}

\title[Statistics of the Vorono\"\i\ cell perimeter in large bi-pointed maps]{Statistics of the Vorono\"\i\ cell perimeter in large bi-pointed maps}
\author{Emmanuel Guitter}
\address{Institut de physique th\'eorique, Universit\'e Paris Saclay, CEA, CNRS, F-91191 Gif-sur-Yvette}
\email{emmanuel.guitter@ipht.fr}

\begin{document}
\maketitle

\begin{abstract}
We study the statistics of the Vorono\"\i\ cell perimeter in large bi-pointed planar quadrangulations. Such maps
have two marked vertices at a fixed given distance $2s$ and their Vorono\"\i\ cell perimeter is simply the length of
the frontier which separates vertices closer to one marked vertex than to the other. 
We characterize the statistics of this perimeter as a function of $s$ 
for maps with a large given volume $N$ both in the scaling limit where
$s$ scales as $N^{1/4}$, in which case the Vorono\"\i\ cell perimeter scales as $N^{1/2}$, 
and in the local limit where $s$ remains finite, in which case
the perimeter scales as $s^2$ for large $s$. The obtained laws are universal and are characteristics of the Brownian map 
and the Brownian plane respectively.   
 \end{abstract}

\section{Introduction}
\label{sec:introduction}
The statistics of planar maps, which are graphs embedded in the sphere considered up to continuous deformations, 
has raised much attention over the years and its study gave rise to a myriad of combinatorial results which describe these mathematical
objects at the discrete level, as well as of probabilistic results which characterize their universal continuous limits such as, depending on the underlying scaling regime, the Brownian map or the Brownian plane. Among refined descriptions of planar maps is the characterization of their 
\emph{Vorono\"\i\ cells} which, for maps with a number of distinguished marked vertices, are the domains which, using the natural graph distance,
regroup vertices according to which marked vertex they are closer to. The simplest realization of Vorono\"\i\ cells is for \emph{bi-pointed maps},
i.e. maps with exactly two marked distinct and distinguished vertices, say $v_1$ and $v_2$, in which case the map is naturally split into two
connected complementary domains gathering vertices closer to $v_1$ than to $v_2$ and conversely.  A natural question is then to characterize 
the statistics of these two Vorono\"\i\ cells in an ensemble where, say the total volume (equal for instance to the number $N$ of faces)
of the map is fixed. A first result was obtained in \cite{G17a} where the law for the \emph{proportion $\phi$ of the total volume $N$}  ($0\leq\phi\leq 1$)
spanned by one of the cells was obtained exactly in the limit of large $N$. It was found that this law is \emph{uniform} between $0$ and $1$,
proving a particular instance of some more general conjecture by Chapuy for Vorono\"\i\ cells within maps of fixed genus \cite{GC16}.
Another result concerned the characterization of Vorono\"\i\ cells for bi-pointed maps where the graph distance 
$d(v_1,v_2)$ is also fixed \cite{G17b}.
In the limit where $N\to \infty$ and $d(v_1,v_2)$ is kept finite, it was shown that exactly one of the two cells has an infinite volume, while 
the other remains finite. The universal law for this finite Vorono\"\i\ cell volume in the limit of large $d(v_1,v_2)$ was derived in \cite{G17b}.

In this paper, we address another question about bi-pointed maps  by looking at their \emph{Vorono\"\i\ cell perimeter}, which is 
the length $L$ of the frontier which separates the two Vorono\"\i\ cells. We again consider the ensemble of bi-pointed maps with a fixed volume 
$N$ and a fixed distance $d(v_1,v_2)$. Explicit  universal expressions are obtained for the Vorono\"\i\ cell perimeter statistics in the 
following
two regimes for which $N\to\infty$: (1) the so-called scaling limit where
$d(v_1,v_2)$ scales as $N^{1/4}$, and (2) the so called local limit where $d(v_1,v_2)$ remains finite.
In case (1), whose behavior characterizes the so-called Brownian map, the Vorono\"\i\ cell perimeter scales as $N^{1/2}$ and
an explicit universal law is derived for the properly rescaled perimeter $L/N^{1/2}$. In case (2), whose behavior
for large $d(v_1,v_2)$ characterizes the so-called Brownian plane, a universal non-trivial law is obtained 
if we now rescale the perimeter $L$ by $d(v_1,v_2)^2$.

The paper is organized as follows: we first discuss in Section 2 the maps under study and give a precise definition of their
Vorono\"\i\ cell perimeter, whose statistics is encoded in some appropriate generating function (Section 2.1). To this generating function
is naturally associated some
\emph{scaling function}, introduced in Section 2.2, whose determination is the key of our study as it may be used to encode the
$N\to \infty$ limit of the Vorono\"\i\ cell perimeter statistics. Section 3 is devoted to the 
computation of this scaling function thanks to a series of differential equations (Section 3.1) solved successively in Sections 3.2 and
3.3. The net result is an explicit expression for the desired scaling function which is then used in Sections 4 and 5 to 
address the Vorono\"\i\ cell perimeter statistics in detail. Section 4 deals with the scaling limit (1) introduced above. 
Section 4.1 explains in all generality how to obtain expected values characterizing the Vorono\"\i\ cell perimeter in this regime  directly from the
scaling function while Section 4.2 presents the corresponding explicit formulas and plots.
Section 5 deals with the local limit (2): here again, we first explain in Section 5.1 the general formalism which 
connects the local limit statistics to the already computed scaling function. This is then used in Section 5.2 to provide
a number of laws, which we make explicit and/or plot for illustration. We present our conclusions in Section 6 while
some technical intermediate computations are discussed in Appendix A. 
 

\section{Bi-pointed maps and Vorono\"\i\ cell perimeter}
\label{sec:voronoi}
\subsection{Generating functions}
This paper deals with the statistics of \emph{bi-pointed planar quadrangulations}, which are planar maps whose all faces have degree $4$, with two marked distinct vertices, distinguished as $v_1$ and $v_2$, at some \emph{even} graph distance from each other. 
Such maps are naturally split into two Vorono\"\i\ cells, which are complementary connected domains
gathering vertices which are closer to one marked vertex than to the other. A precise definition of the Vorono\"\i\ cells is given 
by the Miermont bijection \cite{Miermont2009} which provides a one-to-one correspondence between the above defined bi-pointed planar quadrangulations and
 so-called planar \emph{iso-labelled two-face maps (i-l.2.f.m)}. Those are planar maps with exactly two faces, distinguished as 
$f_1$ and $f_2$, and whose vertices carry positive integer labels. The labels are required to satisfy the following three constraints:
\begin{enumerate}[$(\hbox{L}_1)$]
\item{labels on adjacent vertices differ by $0$ or $\pm1$;}
\item{the minimum label for the set of vertices incident to $f_1$ is $1$;}
\item{the minimum label for the set of vertices incident to $f_2$ is $1$.}
\end{enumerate}
\begin{figure}
\begin{center}
\includegraphics[width=10cm]{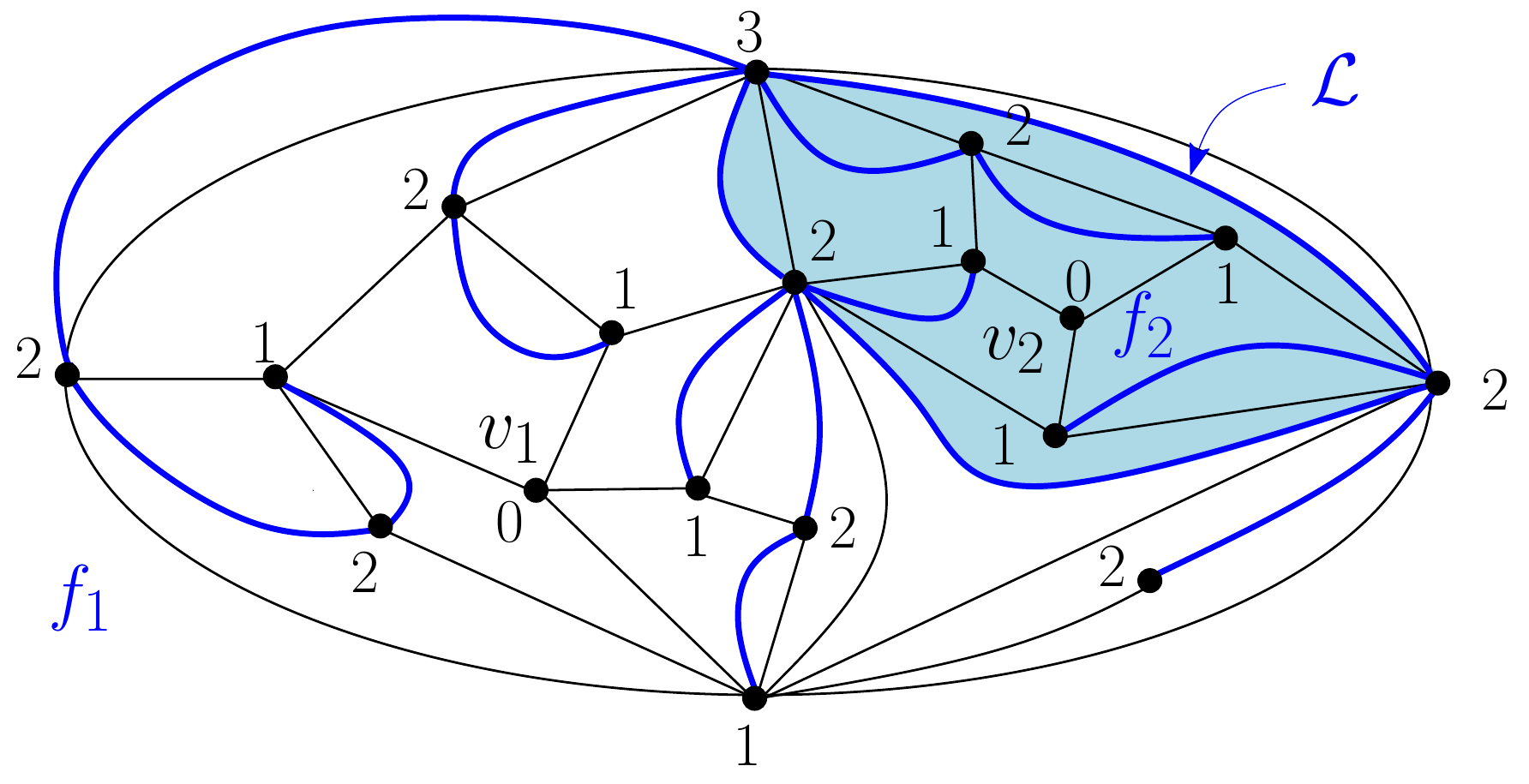}
\end{center}
\caption{An example of planar bi-pointed quadrangulation (thin black edges) whose two marked vertices $v_1$ and $v_2$ are at even distance $2s$ (here $s=2$). Each vertex is labelled by its distance to the closest marked vertex $v_1$ or $v_2$.
Applying the Miermont bijection transforms the quadrangulation into an i-l.2.f.m (thick blue edges) whose two faces $f_1$ and $f_2$ define two Vorono\"\i\ cells (here
the Vorono\"\i\ cell corresponding to $f_2$ has been filled with light blue). The Vorono\"\i\ cell perimeter $L$ is the length of the loop $\mathcal{L}$ separating the 
two faces of the i-l.2.f.m (i.e.\ the boundary length of the light blue domain, here $L=3$).}
\label{fig:cells}
\end{figure}
As it appears in the Miermont bijection, the vertices of the i-l.2.f.m are in one-to-one correspondence with the unmarked (i.e. other than $v_1$
and $v_2$) vertices in the associated bi-pointed quadrangulation and their label corresponds precisely to their graph distance in the quadrangulation 
to the \emph{closest} marked vertex $v_1$ or $v_2$. Moreover, vertices incident to $f_1$ (resp.\ $f_2$) in the i-l.2.f.m correspond to vertices closer to $v_1$ than to $v_2$ (resp.\ closer to $v_2$ than to $v_1$) in the quadrangulation. As a consequence, we define the Vorono\"\i\ cells of the bi-pointed quadrangulation as the two domains spanned by $f_1$ and $f_2$ respectively upon drawing the i-l.2.f.m and the associated quadrangulation on top of each other (see figure \ref{fig:cells}). Each face of the quadrangulation is actually traversed by exactly one edge of the i-l.2.f.m. These edges form a 
simple closed loop $\mathcal{L}$ separating the two faces $f_1$ and $f_2$, together with a number of subtrees attached to the 
vertices of $\mathcal{L}$, on both sides of the loop. Subtree edges of either side correspond to faces in the bi-pointed quadrangulation
which lie within a single Vorono\"\i\ cell, while loop edges correspond instead to faces covered by both cells. 
The length $L$ ($=$number of edges) of the loop $\mathcal{L}$ is called the Vorono\"\i\ cell \emph{perimeter} as it is a natural measure of the size of the boundary between the two cells.

In this paper, we shall study the statistics of $L$ for maps with a \emph{fixed value} $d(v_1,v_2)=2s$ ($s\geq 1$) \emph{of the graph distance between
the two marked vertices $v_1$ and $v_2$.}
In the associated i-l.2.f.m language, this distance is fixed via the following fourth requirement:
\begin{enumerate}[$(\hbox{L}_4)$]
\item{the minimum label for the set of vertices incident to $\mathcal{L}$ is $s$.}
\end{enumerate}

For a fixed graph distance $d(v_1,v_2)=2s$, we may keep a control on both the volume $N$ ($=$ number of faces) of the planar bi-pointed quadrangulations and on their Vorono\"\i\ cell perimeter $L$ by considering their \emph{generating function} $F_s(g,z)$ with a weight 
$g^{N}\, z^{L}$. Via the Miermont bijection, $F_s(g,z)$ is also the generating function of i-l.2.f.m satisfying $(\hbox{L}_4)$ with a weight 
$g^{E}\, z^{L}$, where $E$ is the number of edges of the i-l.2.f.m and $L$ the length of the loop $\mathcal{L}$.

\begin{figure}
\begin{center}
\includegraphics[width=10cm]{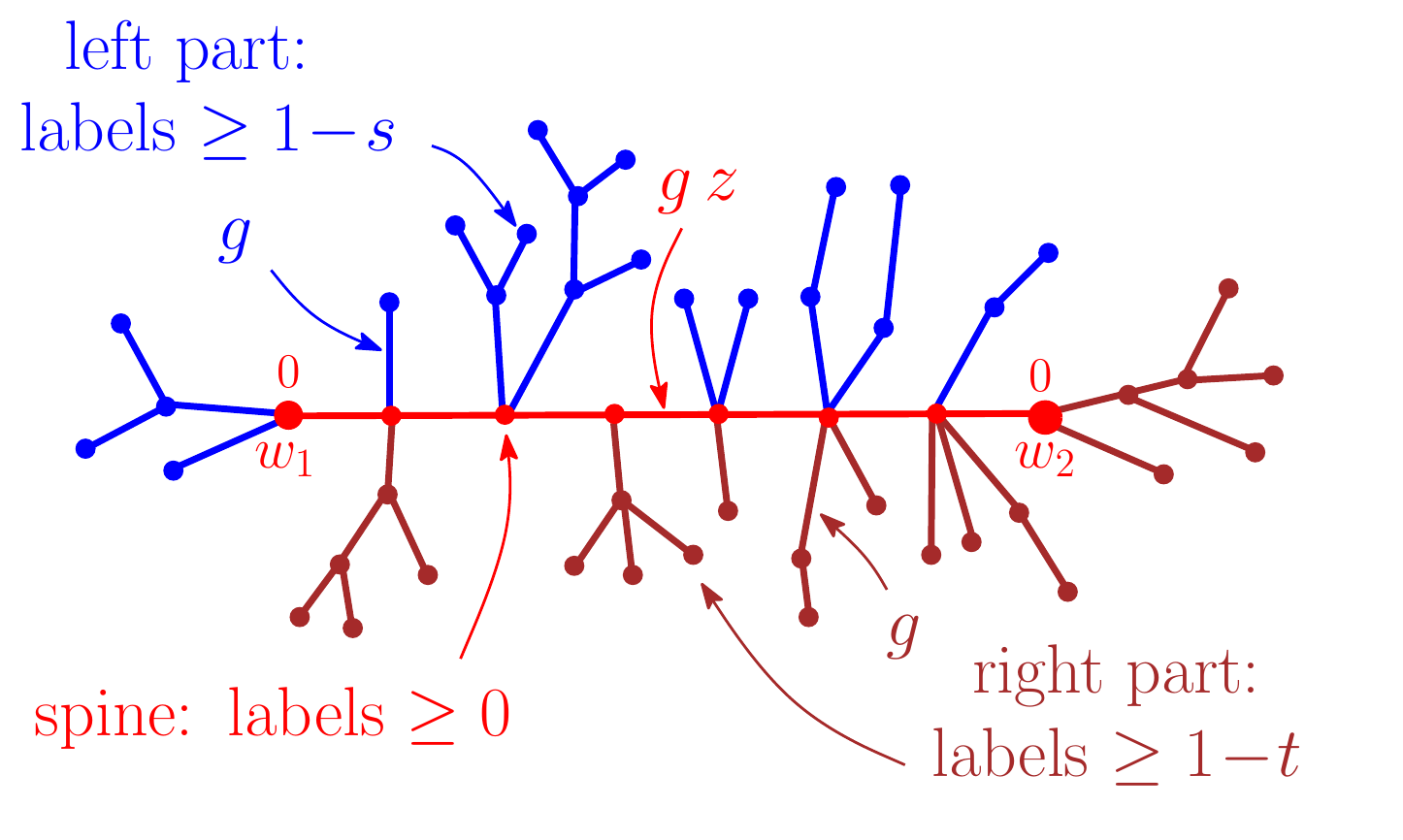}
\end{center}
\caption{A schematic picture of a l.c. enumerated by $X_{s,t}(g,z)$, see text.}
\label{fig:lc}
\end{figure}
As explained in \cite{G17a} in a slightly different context, the generating function $F_s(g,z)$ is closely related to some other generating function, 
$X_{s,t}(g,z)$, of properly weighted \emph{labelled chains (l.c)}, which are planar one-face maps whose vertices carry integer labels and with 
two distinct (and distinguished) marked vertices $w_1$ and $w_2$. The unique shortest path in the map going from $w_1$ to $w_2$ 
forms what is called the \emph{spine} of the map. Apart from its spine, the map consists of (possibly empty) subtrees attached to the left and to
the right of each internal spine vertex and at the marked vertices $w_1$ and $w_2$ (see figure \ref{fig:lc}). 
Given two positive integers $s$ and $t$, $X_{s,t}\equiv X_{s,t}(g,z)$ is the generating function of planar l.c subject to the following constraints:
 
\begin{enumerate}[$(\hbox{$\ell$}_1)$]
\item{labels on adjacent vertices differ by $0$ or $\pm1$;}
\item{$w_1$ and $w_2$ have label $0$. The minimum label for spine vertices is $0$. Spine edges receive a weight $g\, z$;}
\item{the minimum label for vertices of the \emph{left part}, formed by the subtree attached to $w_1$ and all the subtrees attached to the left of any internal spine vertex (with the spine oriented from $w_1$ to $w_2$), is \emph{larger than or equal to} $1-s$. The edges of all these subtrees receive a weight $g$;}
\item{the minimum label for vertices of the \emph{right part}, formed by the subtree attached to $w_2$ and all the subtrees attached to the right of any internal spine vertex, is  \emph{larger than or equal to} $1-t$. The edges of all these subtrees receive a weight $g$;}
\end{enumerate}
For convenience, $X_{s,t}$ includes a first additional term $1$ and, for $s,t>0$, we set $X_{s,0}=X_{0,t}=X_{0,0}=1$.

By some appropriate decomposition of the i-l.2.f.m, it is then easy to relate $F_s(g,z)$ to $X_{s,t}(g,z)$ (see \cite{G17a}
for a detailed description of the decomposition) via:
\begin{equation}
F_s(g,z)= \big(\Delta_s\Delta_t \log(X_{s,t}(g,z))\big)\Big\vert_{t=s} =\log\left(\frac{X_{s,s}(g,z)X_{s-1,s-1}(g,z)}{X_{s-1,s}(g,z)X_{s,s-1}(g,z)}\right)\ ,
\label{eq:FsX}
\end{equation}
where $\Delta_s$ denotes the finite difference operator ($\Delta_s f(s) \equiv f(s)-f(s-1)$). In brief, the above relation is obtained by
first shifting all the labels of the i-l.2.f.m by $-s$ so that the minimum label on the loop $\mathcal{L}$ becomes $0$ and the minimum label
within $f_1$ becomes $1-s$, as well as the minimum label within $f_2$. One then releases the constraints by demanding that the minimum label
within $f_1$ be $\geq 1-s$ and that within $f_2$ be $\geq 1-t$ for some $t$ not necessarily equal to $s$. The loop of the obtained two-face map
forms (by cutting it at each loop vertex labelled $0$) a \emph{cyclic sequence} of particular l.c with no internal spine vertex labelled $0$. Since 
linear (i.e.\ non-cyclic) sequences of such
l.c are precisely enumerated by $X_{s,t}(g,z)$, the desired cyclic sequences are easily seen to be enumerated by $\log(X_{s,t}(g,z))$.
Applying the finite difference operator $\Delta_s$ restores the constraint that the minimum label
within $f_1$ is exactly $1-s$, and similarly,  the finite difference operator $\Delta_t$ restores the constraint that the minimum label
within $f_2$ is exactly $1-t$, which eventually equals the desired value $1-s$ by setting $t=s$.

It is therefore sufficient to obtain an expression for $X_{s,t}(g,z)$ to determine the desired $F_s(g,z)$.
The generating function $X_{s,t}(g,z)$ itself is fully determined by its relation with the generating function $R_s(g)$ for so called well-abelled 
trees\footnote{To be precise, $R_s(g)$ enumerates planted trees with vertices labeled by integers satisfying ($\ell_1$), with root vertex labelled by $0$ and with all labels $\geq 1-s$. Each edge of the tree gets a weight $g$.}(accounting for the subtrees attached to the spine), namely (see \cite{G17a}):
\begin{equation}
X_{s,t}(g,z)=1+g\, z\, R_s(g)R_t(g)X_{s,t}(g,z)\left(1+g\, z\, R_{s+1}(g)R_{t+1}(g)X_{s+1,t+1}(g,z)\right)
\label{eq:eqforXst}
\end{equation}
for $s,t\geq 0$. As for the generating function $R_s(g)$ itself (which depends on $g$ only), its explicit expression is well-known and
takes the parametric form \cite{GEOD}: 
\begin{equation}
R_s(g)=\frac{1+4x+x^2}{1+x+x^2}\frac{(1-x^s)(1-x^{s+3})}{(1-x^{s+1})(1-x^{s+2})}
\quad \hbox{for}\ g= x\frac{1+x+x^2}{(1+4x+x^2)^2}\ .
\label{eq:eqforRs}
\end{equation}
Here $x$ varies in the range $0\leq x\leq 1$ so that $g$ varies in the range $0\leq g\leq 1/12$ (ensuring a proper convergence of the 
generating function at hand). From the explicit form \eqref{eq:eqforRs} and the equations \eqref{eq:eqforXst} and \eqref{eq:FsX}, we may
in principle determine $X_{s,t}(g,z)$ and $F_s(g,z)$. For instance, 
for $z=1$, the solution of \eqref{eq:eqforXst} reads
\begin{equation}
X_{s,t}(g,1)=\frac{(1-x^3)(1-x^{s+1})(1-x^{t+1})(1-x^{s+t+3})}{(1-x)(1-x^{s+3})(1-x^{t+3})(1-x^{s+t+1})}
\label{eq:exactXgone}
\end{equation}
so that
\begin{equation}
F_{s}(g,1)= \log \left\{\frac{\left(1-x^{2 s+3}\right) \left(1-x^{2 s}\right)^2 }{\left(1-x^{2 s-1}\right) \left(1-x^{2 s+2}\right)^2}\right\} \ .
\label{eq:exactFgone}
\end{equation}
For arbitrary $z$ however, we have not been able to obtain such explicit expressions for $X_{s,t}(g,z)$ or $F_s(g,z)$ and we therefore 
have no prediction on the statistics of Vorono\"\i\ cell perimeters for a finite volume $N$ of the quadrangulation. Fortunately, as explained in details in the next section, all the above generating functions have appropriate scaling limits whose knowledge will allow us to characterize the statistics 
of Vorono\"\i\ cell perimeters in quadrangulations \emph{with a large or infinite volume $N$}. 

Before turning to scaling functions, let us end this section by 
looking, now for arbitrary $z$, at the limit 
\begin{equation*}
X_{\infty,\infty}(g,z)\equiv \lim_{s,t\to \infty} X_{s,t}(g,z)\ .
\end{equation*}
It is easily obtained from \eqref{eq:eqforXst} as the solution of 
\begin{equation*}
X_{\infty,\infty}(g,z)=1+z\, \frac{x}{1+x+x^2}\, X_{\infty,\infty}(g,z)\left(1+z\, \frac{x}{1+x+x^2}\, X_{\infty,\infty}(g,z)\right)\ ,
\end{equation*}
which tends to $1$ when $g\to 0$, namely
 \begin{equation}
X_{\infty,\infty}(g,z)= \frac{(1\!+\!x\!+\!x^2) \left(1\!+\!x\!+\!x^2\!-\!z\, x\!-\!\sqrt{(1\!+\!x\!+\!x^2\!-\!z\, x)^2-4\,x^2\,z^2}\right)}{2\,x^2\,z^2}
\label{eq:Xinfinity}
\end{equation}
with the same parametrization of $g$ as in \eqref{eq:eqforRs}. Note that $X_{\infty,\infty}(g,z)$ is well defined for $0\leq z\leq (1+x+x^2)/(3x)$,
hence for $0\leq z\leq 1$ in the whole range $0\leq g\leq 1/12$ ($0\leq x\leq 1$).

\subsection{Scaling functions}
The generating functions $R_s(g)$, $X_{s,t}(g,z)$ and $F_s(g,z)$ are singular when $g$ and $z$ tend simultaneously
toward their critical value $1/12$ and $1$ respectively. The corresponding singularities are encoded in \emph{scaling functions},
as defined below, whose knowledge will allow us to describe the statistics of Vorono\"\i\ cell perimeters in quadrangulations 
whose volume $N$ is large. 
A non trivial scaling regime is obtained by letting $g$ and $z$ approach their critical values as
\begin{equation}
g=G(a,\epsilon)\equiv \frac{1}{12}\left(1-\frac{a^4}{36}\epsilon^4\right)\ ,\quad z=Z(\lambda,\epsilon)\equiv 1-\lambda\, \epsilon^2
\label{eq:gzscal}
\end{equation}
with $\epsilon \to 0$, and letting simultaneously $s$ and $t$ become large as
\begin{equation*}
s=\frac{S}{\epsilon}\ , \qquad t=\frac{T}{\epsilon}\ ,
\end{equation*}
with $a$, $\lambda$, $S$ and $T$ finite positive reals. 
From \eqref{eq:eqforRs}, we have for instance the expansion 
\begin{equation}
R_{\left\lfloor S/\epsilon\right\rfloor}(G(a,\epsilon)) =2+r(S,a)\ \epsilon^2+O(\epsilon^3)\ , \qquad \ r(S,a)=-\frac{a^2 \left(1+10 e^{-a S}+e^{-2 a S}\right)}{3 \left(1-e^{-a S}\right)^2}\ .
 \end{equation}
Expanding \eqref{eq:eqforXst} at second order in $\epsilon$, we have similarly
 \begin{equation*}
X_{\left\lfloor S/\epsilon\right\rfloor,\left\lfloor T/\epsilon\right\rfloor}(G(a,\epsilon),Z(\lambda,\epsilon))  =3+x(S,T,a,\lambda)\ \epsilon+O(\epsilon^2)\ ,
\end{equation*}
where the scaling function $x(S,T,a,\lambda)$ is solution of the non-linear partial
differential equation 
\begin{equation}
2 \big(x(S,T,a,\lambda)\big)^2+6 \big(\partial_Sx(S,T,a,\lambda)+\partial_Tx(S,T,a,\lambda)\big)+27 \big(r(S,a)+r(T,a)\big)=54 \lambda
\label{eq:eqforxstab}
\end{equation}
(with the short-hand notations $\partial_S\equiv \partial/\partial S$ and $\partial_T\equiv \partial/\partial T$). For instance, from the explicitly expression \eqref{eq:exactXgone} at $z=1$ (i.e.\ $\lambda=0$), we immediately deduce
\begin{equation*}
x(S,T,a,0) =-3\, a-\frac{6 a \left(e^{-a S}+e^{-a T}-3 e^{-a (S+T)}+e^{-2a (S+T)}\right)}{\left(1-e^{-a S}\right) \left(1-e^{-a T}\right) \left(1-e^{-a (S+T)}\right)}\ ,
\end{equation*}  
which, as easily checked, satisfies \eqref{eq:eqforxstab} when $\lambda=0$.

Finally, eq.~\eqref{eq:FsX} implies that
\begin{equation}
\begin{split}
F_{\left\lfloor S/\epsilon\right\rfloor}\left(G(a,\epsilon),Z(\lambda,\epsilon)\right) &  = \mathcal{F}(S,a,\lambda)\, \epsilon^3+O(\epsilon^4)\\
\hbox{where} & \quad \mathcal{F}(S,a,\lambda) = \frac{1}{3}\ \partial_S\partial_T x(S,T,a,\lambda)\Big\vert_{T=S}\ .\\
\end{split}
\label{eq:xtoF}
\end{equation}
The next section is devoted to the actual computation of the scaling function $\mathcal{F}(S,a,\lambda)$. We will then show in the following sections
how to use this result to address the question of the Vorono\"\i\ cell perimeter statistics in large maps. 

\section{Computation of the scaling function $\mathcal{F}(S,a,\lambda)$} 

\subsection{Differential equations}
We have not been able to solve the equation \eqref{eq:eqforxstab}, as such, when $\lambda\neq 0$ and therefore have
no explicit expression for $x(S,T,a,\lambda)$ for arbitrary $\lambda$ and for independent $S$ and $T$. Such a formula is however not necessary to compute 
$\mathcal{F}(S,a,\lambda)$ as we may instead recourse to a sequence of two ordinary differential equations, both inherited 
from \eqref{eq:eqforxstab}, which determine successively the following two quantities:
\begin{equation}
\begin{split}
& Y(S,a,\lambda)\equiv x(S,T,a,\lambda)\vert_{T=S}\ ,\\
& H(S,a,\lambda)\equiv \partial_S \partial_T x(S,T,a,\lambda)\vert_{T=S}\ .\\
\end{split}
\label{eq:YJH}
\end{equation}
Note that these quantities depend on a single ``distance parameter'' $S$ instead of two (since we always eventually fix $T=S$). Note also
that the knowledge of $H(S,a,\lambda)$ is exactly what we look for since, from \eqref{eq:xtoF}, we have
\begin{equation*}
\mathcal{F}(S,a,\lambda)=\frac{1}{3}\, H(S,a,\lambda)\ .
\end{equation*}
Since $\frac{d}{dS}Y(S,a,\lambda)=\big(\partial_Sx(S,T,a,\lambda)+\partial_Tx(S,T,a,\lambda)\big)\vert_{T=S}$, we immediately deduce from \eqref{eq:eqforxstab} the following differential equation for $Y(S,a,\lambda)$:
\begin{equation}
2 \big(Y(S,a,\lambda)\big)^2+6 \frac{d}{dS}Y(S,a,\lambda)+54 r(S,a)=54 \lambda\ .
\label{eq:eqforYsab}
\end{equation}
This equation, with appropriate boundary conditions, will allow us to compute $Y(S,a,\lambda)$ explicitly.

Similarly,  differentiating \eqref{eq:eqforxstab} with respect to $S$ and $T$ and then setting $T=S$, we deduce
\begin{equation}
4 Y(S,a,\lambda)\, H(S,a,\lambda)+ \left(\frac{d}{dS}Y(S,a,\lambda)\right)^2+6 \frac{d}{dS}H(S,a,\lambda)=0 \ ,
\label{eq:eqforHsab}
\end{equation}
where we used the identity $\partial_S x(S,T,a,\lambda)\vert_{T=S}=\partial_T x(S,T,a,\lambda)\vert_{T=S}=\frac{1}{2} \frac{d}{dS}Y(S,a,\lambda)$,
obvious by symmetry\footnote{We indeed clearly have $X_{s,t}(g,z)=X_{t,s}(g,z)$, hence $x(S,T,a,\lambda)=x(T,S,a,\lambda)$.}. Knowing $Y(S,a,\lambda)$, this equation will determine $H(S,a,\lambda)$ by some appropriate choice of boundary conditions.

Let us now solve the two equations \eqref{eq:eqforYsab} and \eqref{eq:eqforHsab} successively.

\subsection{Solution of \eqref{eq:eqforYsab}}
The non-linear differential equation \eqref{eq:eqforYsab} is of the Riccati type and is therefore easily transformed, by setting
\begin{equation}
Y(S,a,\lambda)=3 \frac{d}{dS} \log\big(V(S,a,\lambda)\big) \ ,
\label{eq:YV}
\end{equation}
into a second order homogeneous linear differential equation 
\begin{equation*}
\frac{d^2}{dS^2} V(S,a,\lambda)+3(r(S,a)-\lambda)V(S,a,\lambda)=0\ .
\end{equation*}
The solution of this equation is easily found to be 
\begin{equation}
V(S,a,\lambda)=\widetilde{V}(\sigma,A)\quad \hbox{with}\ \sigma\equiv e^{-a S}\ \hbox{and}\ A\equiv \sqrt{1+\frac{3 \lambda}{a^2}}
\label{eq:VVtilde}
\end{equation}
where $\widetilde{V}(\sigma,A)$ is solution of
\begin{equation*}
\sigma\, \frac{d}{d\sigma}\widetilde{V}(\sigma,A)+\sigma^2  \frac{d^2}{d\sigma^2}\widetilde{V}(\sigma,A)=\left(A^2+\frac{12 \sigma }{(1-\sigma )^2}\right) \widetilde{V}(\sigma,A)\ ,
\end{equation*}
hence given by
\begin{equation*}
\begin{split}
&\widetilde{V}(\sigma,A)=\frac{V_0}{(1-\sigma)^3}\left(\sigma^A\, W(\sigma,A)-\alpha\, \sigma^{-A}\, W(\sigma,-A)\right)\ ,\\
&W(\sigma,A)= 4 A^3(1-\sigma)^3+12 A^2(1-\sigma)^2(1+\sigma)+A(1-\sigma)(11+38\sigma+11\sigma^2)\\
&\qquad  \qquad \ +3(1+\sigma)(1+8\sigma+\sigma^2)\ .\\
\end{split}
\end{equation*}
Note that, apart from the global arbitrary multiplicative constant $V_0$ which is unimportant as it drops out in \eqref{eq:YV}, this solution involves some unknown constant $\alpha$ 
which needs to be fixed.
To this end, we recall that, for $s$ and $t$ tending to $\infty$, $X_{s,t}(g,z)$ tends to $X_{\infty,\infty}(g,z)$ given by \eqref{eq:Xinfinity}.
In particular, we have the expansion 
\begin{equation*}
X_{\infty,\infty}(G(a,\epsilon),Z(\lambda,\epsilon))=3-3 \sqrt{a^2+3 \lambda}\ \epsilon+ O(\epsilon^2)
\end{equation*}
from which we deduce 
\begin{equation*}
Y(S,a,\lambda) \underset{S \to \infty}{\sim} -3 \sqrt{a^2+3 \lambda} =-3 a\, A
\end{equation*}
for $a>0$. Sending $S\to \infty$ for $a>0$ amounts to letting $\sigma\to 0$. From \eqref{eq:YV}, we must therefore ensure
\begin{equation*}
-3 a\, \frac{\sigma}{\widetilde{V}(\sigma,A)} \frac{d}{d\sigma} \widetilde{V}(\sigma,A) \underset{\sigma \to 0}{\to} -3 a\, A
\quad\Longleftrightarrow \quad \frac{\sigma}{\widetilde{V}(\sigma,A)} \frac{d}{d\sigma} \widetilde{V}(\sigma,A) \underset{\sigma \to 0}{\to}  A\ .
\end{equation*}
Since $W(\sigma,A)$ has a finite limit when $\sigma\to 0$, and since $A>0$, this latter requirement is fulfilled only for $\alpha=0$ 
(in all the other cases, the limit is $-A$ instead).
This fixes the value of the parameter $\alpha$ to be $0$, while we may set $V_0=1$ without loss of generality,  so that we eventually get
the expression
\begin{equation}
\begin{split}
&\widetilde{V}(\sigma,A)=\frac{\sigma^A\, W(\sigma,A)}{(1-\sigma)^3}\ ,\\
&W(\sigma,A)= 4 A^3(1-\sigma)^3+12 A^2(1-\sigma)^2(1+\sigma)+A(1-\sigma)(11+38\sigma+11\sigma^2)\\
&\qquad  \qquad \ +3(1+\sigma)(1+8\sigma+\sigma^2)\\
\end{split}
\label{eq:Vtilde}
\end{equation}
and, from \eqref{eq:YV},
\begin{equation}
\begin{split}
&Y(S,a,\lambda)=a\, \widetilde{Y}(\sigma,A)\quad \hbox{with}\ \sigma\equiv e^{-a S}\ \hbox{and}\ A\equiv \sqrt{1+\frac{3 \lambda}{a^2}}\ ,\\
&\widetilde{Y}(\sigma,A)= -3\, \frac{\sigma}{\widetilde{V}(\sigma,A)} \frac{d}{d\sigma}\widetilde{V}(\sigma,A)= -3\, A - 36\, \frac{\sigma\, U(\sigma,A)}{(1-\sigma)\, W(\sigma,A)}\ ,\\
&U(\sigma,A) =2A^2(1-\sigma)^2 + 5A(1-\sigma)(1+\sigma)+3(1+3\sigma+\sigma^2)\\
\end{split}
\label{eq:Yexpl}
\end{equation}
and $W(A,\sigma)$ as in \eqref{eq:Vtilde}.

\subsection{Solution of \eqref{eq:eqforHsab}}
Given $Y(S,a,\lambda)$, the solution of the linear first order equation \eqref{eq:eqforHsab} is given by the general formula
\begin{equation*}
\hskip -1.2cm H(S,a,\lambda)=e^{-\frac{2}{3}\int_{S_0}^S dT\, Y(T,a,\lambda)}\left\{H(S_0,a,\lambda)-\frac{1}{6}\int_{S_0}^S dT\, \Big(\frac{d}{dT}Y(T,a,\lambda)\Big)^2 e^{\frac{2}{3}\int_{S_0}^T dT'\, Y(T',a,\lambda)}\right\}
\end{equation*}
for some arbitrary $S_0>0$. Here the value of $H(S_0,a,\lambda)$ is some integration constant which \emph{may a priori be
chosen arbitrarily}. As we shall see below, it will be fixed in our case by the requirement that $H(S,a,\lambda)$ has the desired large $S$ behavior.
It is straightforward from \eqref{eq:YV} that
\begin{equation*}
e^{-\frac{2}{3}\int_{S_0}^S dT\, Y(T,a,\lambda)}=\left(\frac{V(S_0,a,\lambda)}{V(S,a,\lambda)}\right)^2
\end{equation*}
with $V(S,a,\lambda)$ as in eqs.~\eqref{eq:VVtilde} and \eqref{eq:Vtilde}. 
We may therefore write the solution of \eqref{eq:eqforHsab} as
\begin{equation*}
\hskip -1.2cm H(S,a,\lambda)=\frac{1}{\left(V(S,a,\lambda)\right)^2}\left\{\!\left(V(S_0,a,\lambda)\right)^2\, H(S_0,a,\lambda)-\frac{1}{6}\int_{S_0}^S \!\!\!\!dT\, \Bigg(\!V(T,a,\lambda)\, \frac{d}{dT}Y(T,a,\lambda)\!\Bigg)^2\!\right\}
\end{equation*}
which, upon changing again variables by setting
\begin{equation*}
H(S,a,\lambda)=a^3\, \widetilde{H}(\sigma,A)\quad \hbox{with}\ \sigma\equiv e^{-a S}\ \hbox{and}\ A\equiv \sqrt{1+\frac{3 \lambda}{a^2}}
\end{equation*}
and upon using eqs.~\eqref{eq:Vtilde} and \eqref{eq:Yexpl}, rewrites explicitly 
\begin{equation*}
\begin{split}
&\hskip -1.2cm\widetilde{H}(\sigma,A)=\frac{(1-\sigma)^6}{\sigma^{2A} \left(
W(\sigma,A)\right)^2}\left\{\frac{\sigma_0^{2A}\, \left(W(\sigma_0,A)\right)^2}{(1-\sigma_0)^6}\, \widetilde{H}(\sigma_0,A)+\int_{\sigma_0}^
\sigma d\tau\, M(\tau,A)\right\}\ ,\\
&\hskip -1.2cm M(\sigma,A)= \frac{216\, \sigma^{2 A+1}\left(P(\sigma,A)\right)^2}{(1-\sigma)^{10}\left(W(\sigma,A)\right)^2} \ ,
\\
&\hskip -1.2cm P(\sigma,A)= \Big(W(\sigma,A)(1-\sigma)\Big)^2\, \frac{d}{d\sigma} \left(\frac{\sigma\, U(\sigma,A)}{(1-\sigma)\, W(\sigma,A)}\right)\\
&= 8A^5(1-\sigma)^5 (1+\sigma)+4 A^4(1-\sigma)^4(11+20 \sigma+11\sigma^2)\\
& +2 A^3(1-\sigma)^3 (1+\sigma) (47
 + 86 \sigma +47\sigma^2)\\
 &+A^2 (1-\sigma)^2 (97+386 \sigma
 +594 \sigma^2+386 \sigma^3 +97\sigma^4)\\
 &+6 A(1-\sigma) (1+\sigma) (8+33 \sigma+68
   \sigma^2+33 \sigma^3 +8\sigma^4)\\
   &+9\, (1+6 \sigma+27 \sigma^2+32 \sigma^3+27 \sigma^4+6
   \sigma^5 +\sigma^6)\ , \\
\end{split}
\end{equation*}
where $\sigma_0= e^{-a S_0}$. Remarkably enough, the function $M(\sigma,A)$ above
has an explicit primitive (in $\sigma$) given by:
\begin{equation*}
\begin{split}
&\hskip -1.2cm C(\sigma,A)=\frac{12}{35}\, \sigma^{2A+2}
\Bigg\{\frac{Q(\sigma,A)}{(1-\sigma)^9W(\sigma,A)}\\
& \qquad \qquad \qquad -A(A-1)(A^2-1)(4A^2-1)(4A^2-9)\ _2F_1(1,2A+2,2A+3,\sigma)\Bigg\}\\
\end{split}
\end{equation*}
where $_2F_1(1,2A+2,2A+3,\sigma)$ is a simple hypergeometric function 
\begin{equation}
_2F_1(1,2A+2,2A+3,\sigma)=(2A+2)\, \sum_{n\geq 0} \frac{\sigma^n}{n+2A+2}
\label{eq:def2F1}
\end{equation}
and where $Q(\sigma,A)$ a polynomial of degree $11$ in $\sigma$ given by
\begin{equation}
Q(\sigma,A)=\sum_{i=0}^{11}q_i(\sigma) A^{11-i}(1-\sigma)^{11-i}
\label{eq:Qval}
\end{equation}
with
\begin{equation*}
\begin{split}
&\hskip -1.2cm q_{0}(\sigma)=64\ , \quad q_{1}(\sigma)=32 (4+7 \sigma)\ , \quad q_{2}(\sigma)=80
\, (-3+13 \sigma )\\
&\hskip -1.2cm q_3(\sigma)=48 (-12+3 \sigma +39 \sigma^2-17 \sigma^3)\ , \qquad 
q_4(\sigma)=12  (17-278 \sigma +400 \sigma^2-90 \sigma^3-65 \sigma^4)\\
&\hskip -1.2cm q_5(\sigma)=6 (980+1365 \sigma +424 \sigma^2+1374 \sigma^3-876 \sigma^4+109 \sigma^5)\\
&\hskip -1.2cm q_6(\sigma)= (35363+88857 \sigma +75432 \sigma^2+45698 \sigma^3-1965 \sigma^4-2715 \sigma^5+1250 \sigma^6)\\
&\hskip -1.2cm q_7(\sigma) =(100336+279172 \sigma +377163 \sigma^2+267688
   \sigma^3+118786 \sigma^4-12360 \sigma^5+2611 \sigma^6+244 \sigma^7) \\
&\hskip -1.2cm q_8(\sigma)=2 (74901+220667 \sigma +395067 \sigma^2+388137 \sigma^3+222947 \sigma^4+77061 \sigma^5-2895 \sigma^6\\ & +1355 \sigma^7-240 \sigma^8)\\
   &\hskip -1.2cm q_9(\sigma)= 3  (40659+122280 \sigma +287397
   \sigma^2+350940 \sigma^3+284873 \sigma^4+126322 \sigma^5+37701 \sigma^6\\ & +724 \sigma^7+206  \sigma^8-102 \sigma ^9)\\
   &\hskip -1.2cm q_{10}(\sigma)=9 (5673+16595 \sigma +51926 \sigma^2+83180 \sigma^3+83930
   \sigma^4+50678 \sigma^5+17870 \sigma^6+4994 \sigma^7\\ & +155 \sigma^8+5 \sigma^9-6 \sigma ^{10})\\
   &\hskip -1.2cm q_{11}(\sigma)=27 (315+840 \sigma +3465 \sigma^2+9576 \sigma^3+6510 \sigma^4+9864
   \sigma^5+3150 \sigma^6+1000 \sigma^7\\ & +279 \sigma^8+\sigma^{10})\ .\\
\end{split}
\end{equation*}
Knowing $C(\sigma,A)$, we may thus write
\begin{equation*}
\widetilde{H}(\sigma,A)=\frac{(1-\sigma)^6}{\sigma^{2A} \left(
W(\sigma,A)\right)^2}\left\{\frac{\sigma_0^{2A}\, \left(W(\sigma_0,A)\right)^2}{(1-\sigma_0)^6}\, \widetilde{H}(\sigma_0,A)+C(\sigma,A)-C(\sigma_0,A)\right\}\ ,
\end{equation*}
leaving us with the determination of the free integration constant $\widetilde{H}(\sigma_0,A)=H(S_0,a,\lambda)/a^3$. 
From its combinatorial interpretation, $F_s(g,z)$ is expected to vanish at large $s$ (see for instance \eqref{eq:exactFgone} when $z=1$) 
and so is thus $\mathcal{F}(S,a,\lambda)=(1/3)H(S,a,\lambda)$ at large $S$. Sending $S\to\infty$ amounts (for $a>0$) to let $\sigma\to 0$,
in which case we find
\begin{equation*}
C(\sigma,A)=\sigma^{2A}\Big( 108 \, (A+1)(2A+3)^2\, \sigma^2+O(\sigma^3)\Big)\ .
\end{equation*}
Demanding that $H(S,a,\lambda)\to 0$ at large $S$ forces us to choose the (so far arbitrary) value $\widetilde{H}(\sigma_0,A)$ such that 
$\frac{\sigma_0^{2A}\, \left(W(\sigma_0,A)\right)^2}{(1-\sigma_0)^6}\, \widetilde{H}(\sigma_0,A)-C(\sigma_0,A)=0$ in order to avoid a divergence of 
$\widetilde{H}(\sigma,A)$ as $\sigma^{-2A}$ when $\sigma\to 0$. This leads to the 
final expression\footnote{At small $\sigma$, we then have the expansion 
$\widetilde{H}(\sigma,A)=\frac{108 \, \sigma ^2}{(A+1) (2 A+1)^2}+O\left(\sigma ^3\right)$.}
\begin{equation}
\begin{split}
\hskip -1.2cm \widetilde{H}(\sigma,A)& =\frac{(1-\sigma)^6\, C(\sigma,A)}{\sigma^{2A} \left(
W(\sigma,A)\right)^2}\\
&=\frac{12}{35}\, \sigma^2 
\Bigg\{\frac{Q(\sigma,A)}{(1-\sigma)^3\left(W(\sigma,A)\right)^3}\\
& \qquad   -A(A-1)(A^2-1)(4A^2-1)(4A^2-9)\ \frac{(1-\sigma)^6}{\left(
W(\sigma,A)\right)^2}\ _2F_1(1,2A+2,2A+3,\sigma)\Bigg\}
\end{split}
\label{eq:Hval}
\end{equation}
with $W(A,\sigma)$ as in \eqref{eq:Vtilde} and $Q(\sigma,A)$ as in \eqref{eq:Qval} above. The desired scaling function $\mathcal{F}(S,a,\lambda)$ is then obtained directly via
\begin{equation}
\mathcal{F}(S,a,\lambda)=\frac{1}{3} a^3\, \widetilde{H}(\sigma,A)\ ,  \qquad  \sigma\equiv e^{-a S}\ , \qquad  A\equiv \sqrt{1+\frac{3 \lambda}{a^2}}\ .
\label{eq:Fval}
\end{equation}
\section{Large maps in the scaling limit} 

\subsection{Quadrangulations with a large fixed volume}
\begin{figure}
\begin{center}
\includegraphics[width=13cm]{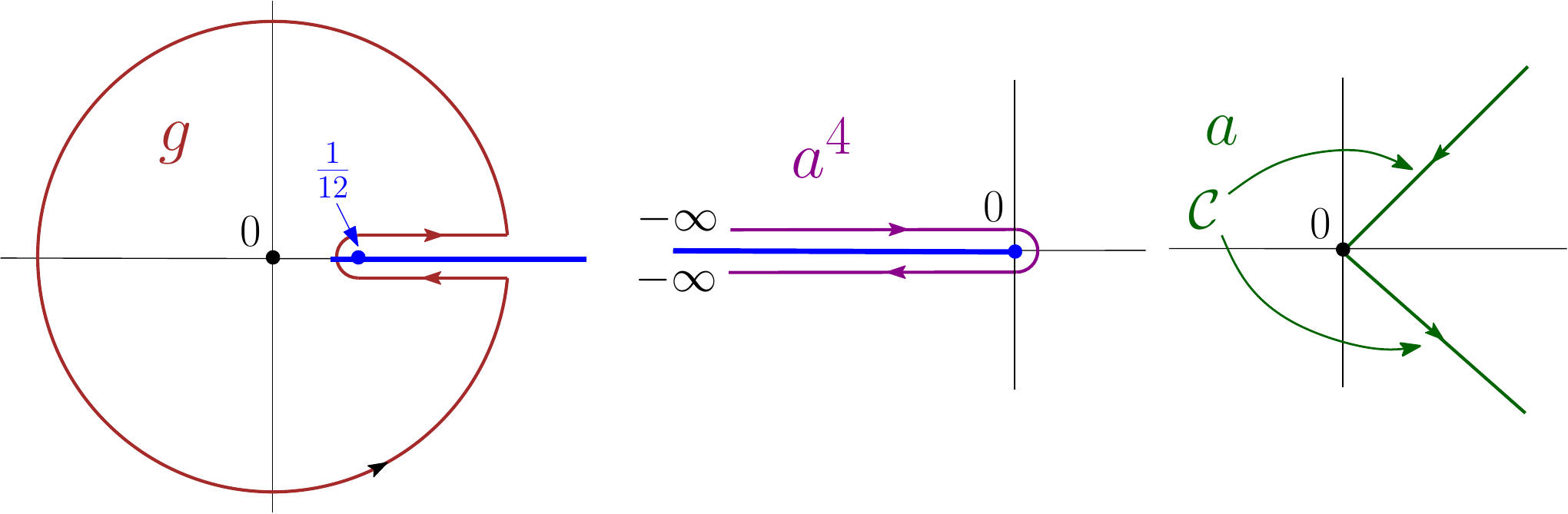}
\end{center}
\caption{The contour integral over $g$ used to extract the coefficient $g^N$ of, say $F_s(g,z)$. In the $g$ plane, the contour around $0$
is deformed so as to surround the cut of $F_s(g,z)$ for real $g>1/12$. By setting $g=G(a,N^{-1/4})$, the excursion
around the cut corresponds in the $a^4$ variable to a path around the negative real axis which makes a back and forth 
excursion from and to $-\infty$ in the limit $N\to\infty$. In the $a$ plane, the corresponding path $\mathcal{C}$ is made of two straight lines at $45^\circ$ from the real axis, oriented as shown from and to infinity. }
\label{fig:contour}
\end{figure}

Fixing the value $N$ of the volume of the quadrangulations amounts to  
extracting the coefficient $g^N$ of the various generating functions. For instance, the \emph{number} of bi-pointed planar quadrangulations of fixed volume $N$ 
and with a fixed distance $2s$ between the marked vertices is directly given by
\begin{equation*}
[g^N]F_s(g,1)\ .
\end{equation*}
As for the statistics of the Vorono\"\i\ cell perimeter $L$ within bi-pointed quadrangulations of fixed volume $N$ 
and fixed distance $2s$ between the marked vertices, it is fully characterized by considering the 
expected value $E_{N,s}(z^L)$ of $z^L$ in this ensemble, for $0\leq z\leq 1$. This expectation is directly given by the simple
ratio
\begin{equation}
E_{N,s}(z^L)=\frac{[g^N]F_s(g,z)}{[g^N]F_s(g,1)}\ .
\label{eq:ensdef}
\end{equation}  
In practice, extracting the coefficient $g^N$ may be performed via a contour integral over $g$ around $0$, which may then
be transformed into an integral over the variable $a$ such that 
\begin{equation*}
g=G\left(a,N^{-1/4}\right)=\frac{1}{12}\left(1-\frac{a^4}{36\, N}\right)\ .
\end{equation*}
The integral over $a$ runs on some contour $\mathcal{C}_N$, which, at large $N$, eventually tends to some simple contour $\mathcal{C}$
made of two straight lines at $45^\circ$ from the real axis, as displayed in figure \ref{fig:contour}.

The \emph{scaling limit} corresponds to let $N$ tend to infinity while looking at distances $2s$ of order $N^{1/4}$. This is achieved  by setting 
$s=S\, N^{1/4}$ where $S$ is kept finite, in which case we have for instance, for $z=1$:
\begin{equation*}
\begin{split}
\hskip -1.2cm [g^N]F_{\left\lfloor S\, N^{1/4} \right\rfloor}(g,1) &=\frac{1}{2{\rm i}\pi}\oint \frac{dg}{g^{N+1}} F_{\left\lfloor S\, N^{1/4} 
\right\rfloor}(g,1)\\ &
\ = \frac{1}{2{\rm i}\pi} \frac{12^N}{N} \int_{\mathcal{C}_N} \frac{-da\, a^3}{9} \frac{1}{\left(1-\frac{a^4}{36\, N}\right)^{N+1}}
F_{\left\lfloor S\, N^{1/4} \right\rfloor}\left(G(a,N^{-1/4}),1\right)\\
& \ \underset{N \to \infty}{\sim}\frac{1}{2{\rm i}\pi} \frac{12^N}{N^{7/4}}  \int_{\mathcal{C}} \frac{-da\, a^3}{9} e^{\frac{a^4}{36}}
\mathcal{F}\left(S,a,0\right)
\end{split}
\end{equation*}
where we have used the expansion \eqref{eq:xtoF} (at $z=1$, i.e.\ at $\lambda=0$) with $\epsilon=1/N^{1/4}$.
Dividing\footnote{In $F_{\left\lfloor S\, N^{1/4} \right\rfloor}(g,1)$, the half-distance is fixed at some precise integer value close to $S\, N^{1/4}$. Letting instead the half-distance vary between $S\, N^{1/4}$ and $(S+dS)\, N^{1/4}$ requires considering the quantity 
$F_{\left\lfloor S\, N^{1/4} \right\rfloor}(g,1)\times N^{1/4}\, dS$.} the quantity $F_{\left\lfloor S\, N^{1/4} \right\rfloor}(g,1)\times N^{1/4}\, dS$ by the large $N$ asymptotic number\footnote{See for instance \cite{G17a} for a derivation of this asymptotic number.} $12^N/(4\sqrt{\pi} N^{3/2})$ of bi-pointed quadrangulations 
having their 
two marked vertices at (arbitrary) even distance from each other leads to the probability $\rho(S) dS$ that two marked vertices
chosen uniformly at random \emph{at some even distance from each other} in a quadrangulation of large volume $N$ have their half-distance between 
$S\, N^{1/4}$ and $(S+dS)\, N^{1/4}$.
The quantity $\rho(S)$ is nothing but the so-called \emph{distance density profile} of the Brownian map, obtained here in the context of bi-pointed
quadrangulations, and is thus given by
\begin{equation}
\begin{split}
\rho(S)&= \frac{2}{{\rm i}\sqrt{\pi}}  \int_{\mathcal{C}} \frac{-da\, a^3}{9} e^{\frac{a^4}{36}}
\mathcal{F}\left(S,a,0\right)\\
&= \frac{2}{{\rm i}\sqrt{\pi}}  \int_0^\infty 2t\, dt\, e^{-t^2} \left\{\mathcal{F}\left(S,(1-{\rm i})\sqrt{3t},0\right)-\mathcal{F}\left(S,(1+{\rm i})\sqrt{3t},0\right)\right\}\\
\end{split}
\label{eq:rho}
\end{equation}
upon changing variable and parametrizing the contour $\mathcal{C}$ by setting $a=(1\mp {\rm i})\sqrt{3t}$, with $t$ real and positive.

We may compute in a similar way the quantity
\begin{equation*}
\begin{split}
&\hskip -1.2cm [g^N]F_{\left\lfloor S\, N^{1/4} \right\rfloor}\left(g,Z(\lambda,N^{-1/4})\right)=\frac{1}{2{\rm i}\pi}\oint \frac{dg}{g^{N+1}} F_{\left\lfloor S\, N^{1/4} 
\right\rfloor}\left(g,Z(\lambda,N^{-1/4})\right)\\ &
\ = \frac{1}{2{\rm i}\pi} \frac{12^N}{N} \int_{\mathcal{C}_N} \frac{-da\, a^3}{9} \frac{1}{\left(1-\frac{a^4}{36\, N}\right)^{N+1}}
F_{\left\lfloor S\, N^{1/4} \right\rfloor}\left(G(a,N^{-1/4}),Z(\lambda,N^{-1/4})\right)\\
& \ \underset{N \to \infty}{\sim} \frac{1}{2{\rm i}\pi} \frac{12^N}{N^{7/4}} \int_{\mathcal{C}} \frac{-da\, a^3}{9} e^{\frac{a^4}{36}}
\mathcal{F}\left(S,a,\lambda\right)\ .
\end{split}
\end{equation*}
Using $Z(\lambda,N^{-1/4})^L\sim e^{-\lambda\, L/N^{1/2}}$ at large $N$, we deduce, after normalization, that
\begin{equation}
\begin{split}
&E_{N,\left\lfloor S\, N^{1/4} \right\rfloor}\left(e^{-\lambda \, L/N^{1/2}}\right) \underset{N \to \infty}{\sim} \frac{\rho(S,\lambda)}{\rho(S)}\ ,
\\
&\rho(S,\lambda)= \frac{2}{{\rm i}\sqrt{\pi}}  \int_0^\infty 2t\, dt\, e^{-t^2} \left\{\mathcal{F}\left(S,(1-{\rm i})\sqrt{3t},\lambda\right)-\mathcal{F}\left(S,(1+{\rm i})\sqrt{3t},\lambda\right)\right\}\\
\end{split}
\label{eq:rholambda}
\end{equation}
and $\rho(S)=\rho(S,0)$ as in \eqref{eq:rho} above. In the scaling limit, a non-trivial law is therefore obtained for the rescaled perimeter $L/N^{1/2}$.
Let us now examine this law in more details.

\subsection{Explicit expressions and plots}
\begin{figure} 
\begin{center}
\includegraphics[width=8cm]{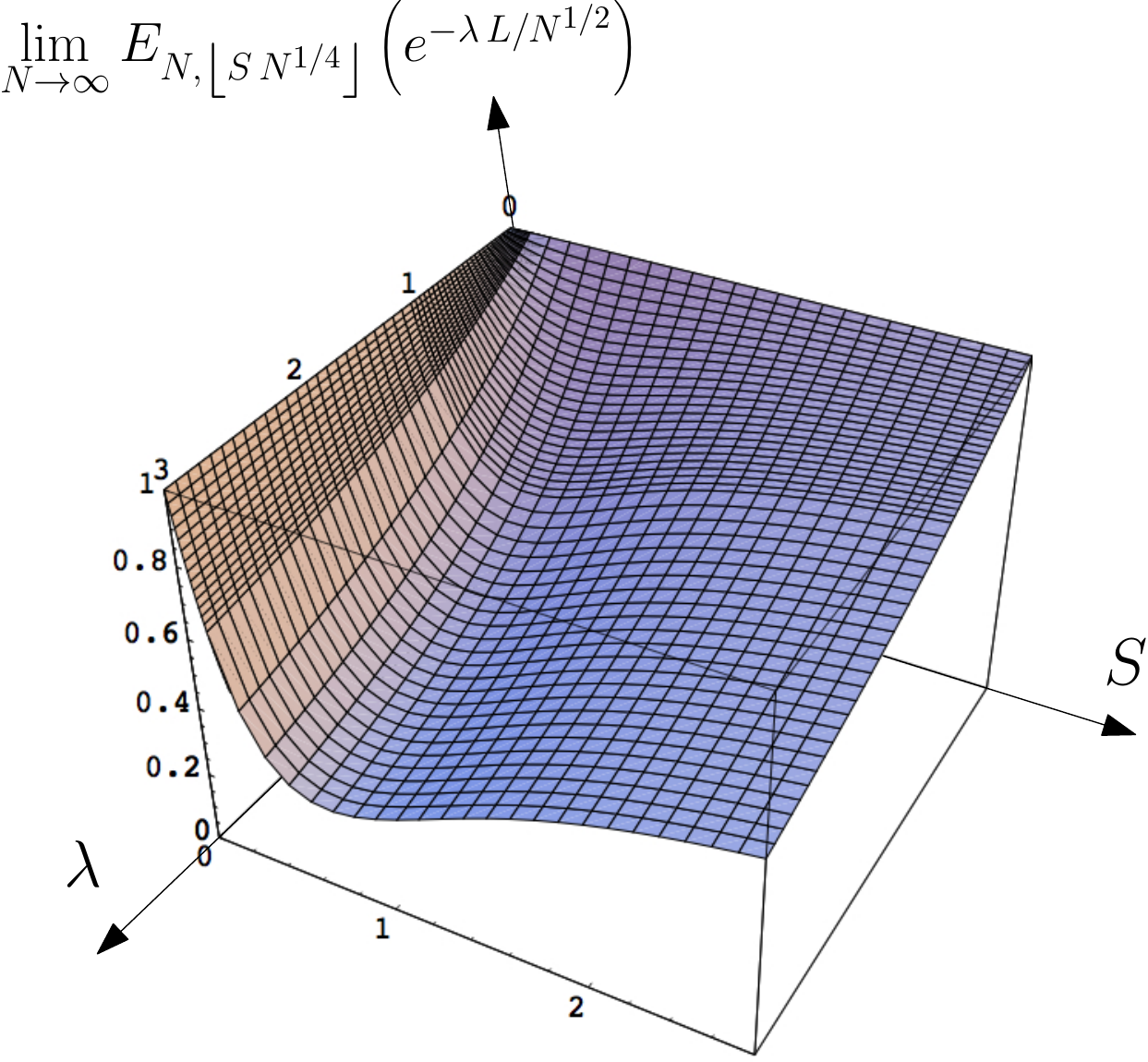}
\end{center}
\caption{A plot of the large $N$ limit of $E_{N,\left\lfloor S\, N^{1/4} \right\rfloor}\left(e^{-\lambda \, L/N^{1/2}}\right)$ as a function of $S$ 
and $\lambda$, as given by \eqref{eq:rholambda}.}
\label{fig:ENSLambda}
\end{figure}
Inserting the explicit expression \eqref{eq:Fval} of $\mathcal{F}(S,a,\lambda)$ into \eqref{eq:rholambda} allows to write
$\rho(S,\lambda)$ as a simple integral over $t$ which may be evaluated numerically. Note that for $S$ and $\lambda$ real positive, 
$\rho(S,\lambda)$ is actually real. Dividing $\rho(S,\lambda)$  by $\rho(S)=\rho(S,0)$ allows us, from \eqref{eq:rholambda}, to then evaluate 
the large $N$ limit of the quantity $E_{N,\left\lfloor S\, N^{1/4} \right\rfloor}\left(e^{-\lambda \, L/N^{1/2}}\right)$ as a function of $S$ 
and $\lambda$. This quantity is plotted for illustration in figure \ref{fig:ENSLambda} and, by definition, lies
between $0$ and $1$. Let us now characterize it in more details via a number of more explicit expressions.
Expanding $\mathcal{F}(S,a,\lambda)$ at first order in $\lambda$, we get
\begin{equation*}
\begin{split}
&\mathcal{F}(S,a,\lambda)=\mathcal{F}(S,a,0)+\lambda\times \frac{\partial}{\partial\lambda}\mathcal{F}(S,a,0)+O(\lambda^2)\ ,\\
&\mathcal{F}(S,a,0)=\frac{2 a^3 \sigma ^2 \left(1+\sigma ^2\right)}{\left(1-\sigma ^2\right)^3}\qquad \hbox{with}\ \sigma=e^{-a\, S}\ ,\\
&\frac{\partial}{\partial\lambda}\mathcal{F}(S,a,0)=-\frac{a \sigma ^2 \left(385-189 \sigma +154 \sigma^2+54 \sigma^3-11 \sigma^4-9 \sigma^5\right)}{70 (1-\sigma ) (1+\sigma )^4} \ .\\
\end{split}
\end{equation*}
From the explicit expression of $\mathcal{F}(S,a,0)$, we may write via \eqref{eq:rho} a rather explicit expression for $\rho(S)$ in the form 
of an integral over $t$.  
Introducing the short hand notations
\def\co{\rm co}
\def\si{\rm si}
\def\ch{\rm ch}
\def\sh{\rm sh}
\begin{equation*}
\ch=\cosh(\sqrt{3t}\, S)\ , \quad
\sh=\sinh(\sqrt{3t}\, S)\ , \quad
\co=\cos(\sqrt{3t}\, S)\ , \quad
\si=\sin(\sqrt{3t}\, S)\ ,
\end{equation*}
we find 
\begin{equation}
\begin{split}
 \hskip -1.2cm \rho(S)=&\sqrt{\frac{3}{\pi }} \int_0^\infty  dt\, 24  t^{5/2}e^{-t^2} \frac{1}{\left(\co^2-\ch^2\right)^3}\times
\\ &\qquad\quad  \Big\{ \ch\, \sh \left(\left(2\co^2\!-\!1\right) \left(\co^2\!+\!\ch^2\!-\!1\right)\!-\!1\right)
\!+\!\co\, \si \left(\left(2\ch^2\!-\!1\right) \left(\co^2\!+\!\ch^2\!-\!1\right)\!-\!1\right) \Big\} .\\
\end{split}
\label{eq:rhoexp}
\end{equation}
\begin{figure}
\begin{center}
\includegraphics[width=13cm]{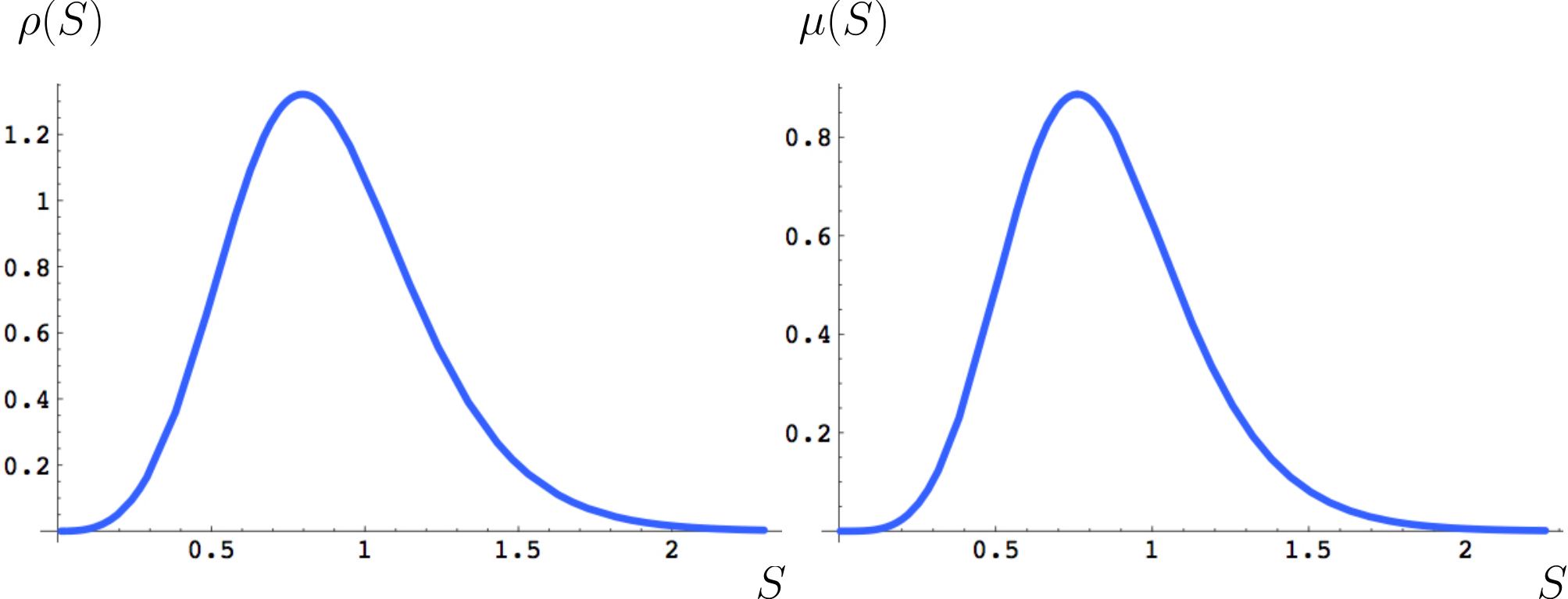}
\end{center}
\caption{Plots of $\rho(S)$ and $\mu(S)$, as given by \eqref{eq:rhoexp} and \eqref{eq:muexp} respectively.}
\label{fig:rhomu}
\end{figure}
This function is plotted in figure \ref{fig:rhomu} and the above formula matches already known expressions in the literature \cite{GEOD} for the distance 
density profile of the Brownian map. We have in particular the small $S$ behavior
\begin{equation*}
\rho(S)=\frac{48}{7} \, S^3+O(S^4)\ .
\end{equation*}
More interesting is the linear term in $\lambda$: expanding eq.~\eqref{eq:rholambda} at first order in $\lambda$ yields indeed
\begin{equation}
\begin{split}
&E_{N,\left\lfloor S\, N^{1/4} \right\rfloor}\left(L/N^{1/2}\right) \underset{N \to \infty}{\sim} \frac{\mu(S)}{\rho(S)}
\\
&\mu(S)= \frac{2}{{\rm i}\sqrt{\pi}}  \int_0^\infty 2t\, dt\, e^{-t^2} \left\{-\frac{\partial}{\partial\lambda}\mathcal{F}\left(S,(1-{\rm i})\sqrt{3t},0\right)+ \frac{\partial}{\partial\lambda}\mathcal{F}\left(S,(1+{\rm i})\sqrt{3t},0\right)\right\}\ .\\
\end{split}
\label{eq:ENexp}
\end{equation}
From the above expression for $\frac{\partial}{\partial\lambda}\mathcal{F}\left(S,a,0\right)$, we arrive at the following expression (with the same short hand notations as before):
\begin{equation}
\begin{split}
&\hskip -1.4cm \mu(S)=\sqrt{\frac{3}{\pi }} \int_0^\infty dt\, \frac{4}{35}  t^{3/2}e^{-t^2} \Bigg\{ 16\, (\co\!-\!\si) (\ch\!-\!\sh)-9\, (\co^2\!-\!2 \co\, \si\!-\!\si^2) (\ch^2\!-\!2 \ch\, \sh\!+\!\sh^2)\\& \qquad \!-\!7
\!-\!12\, \frac{\si\!-\!\sh}{\co\!-\!\ch}\!-\!40\, \frac{\si^2\!-\!2 \si\, \sh\!-\!\sh^2}{(\co\!+\!\ch)^2}
  \!-\!21\, \frac{\si^4\!-\!4 \si^3\, \sh\!-\!6 \si^2\, \sh^2\!+\!4 \si\, \sh^3\!+\!\sh^4}{(\co\!+\!\ch)^4}
\Bigg\}\ .\\
\label{eq:muexp}
\end{split}
\end{equation}
\begin{figure}
\begin{center}
\includegraphics[width=10cm]{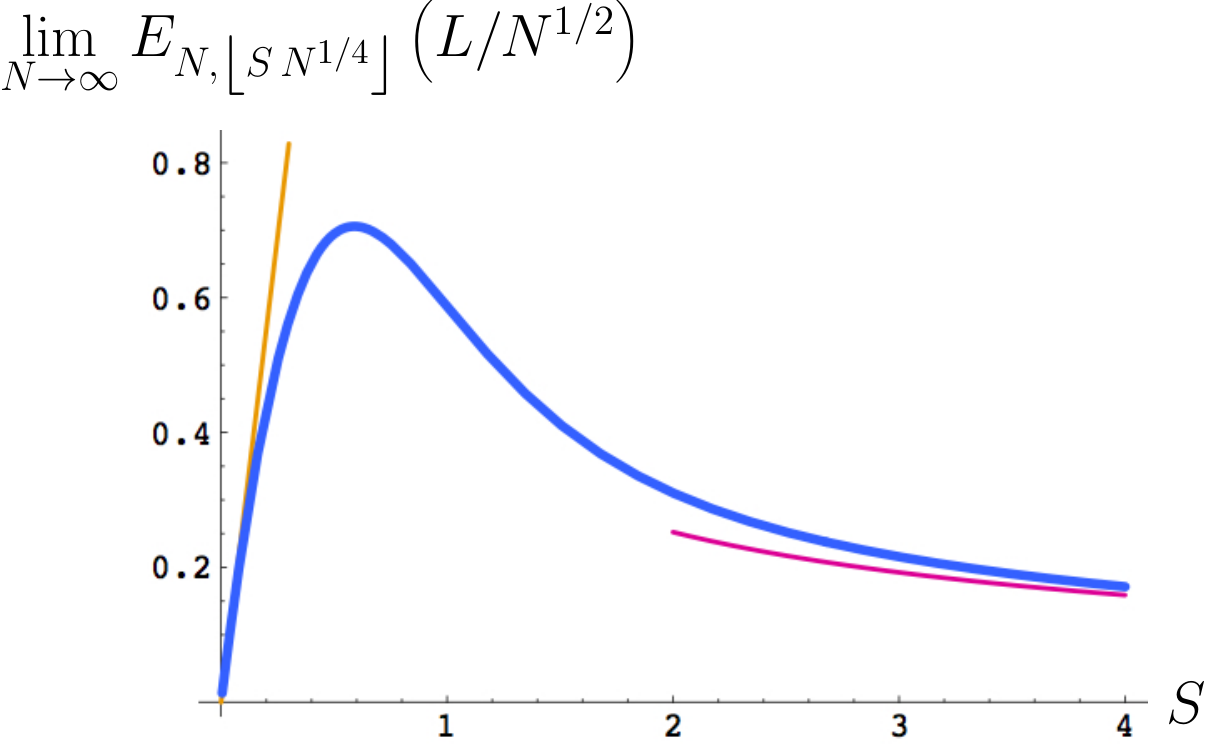}
\end{center}
\caption{A plot of the large $N$ limit of $E_{N,\left\lfloor S\, N^{1/4} \right\rfloor}\left(L/N^{1/2}\right)$, as given by \eqref{eq:rhoexp}--\eqref{eq:ENexp}--\eqref{eq:muexp}. We also indicated its limiting small $S$ linear behavior, as given by \eqref{eq:smallS}, and its large $S$ power law behavior, as given by \eqref{eq:largeS}.}
\label{fig:aveLvsS}
\end{figure}
The function $\mu(S)$ is plotted in figure \ref{fig:rhomu} while the large $N$ expectation 
$E_{N,\left\lfloor S\, N^{1/4} \right\rfloor}\left(L/N^{1/2}\right)=\mu(S)/\rho(S)$ is represented in figure \ref{fig:aveLvsS}. 
At small $S$, we have the behavior
\begin{equation*}
\mu(S)=\frac{171}{10} \sqrt{\frac{3}{\pi }} \Gamma \left(\frac{9}{4}\right) \, S^4+O(S^5)
\end{equation*}
so that
\begin{equation}
E_{N,\left\lfloor S\, N^{1/4} \right\rfloor}\left(L/N^{1/2}\right)\underset{N \to \infty}{\sim} \frac{399}{160} \sqrt{\frac{3}{\pi }} \Gamma \left(\frac{9}{4}\right) \, S+O(S^2)\ .
\label{eq:smallS}
\end{equation}
This limiting behavior is plotted in figure \ref{fig:aveLvsS}. Returning to unscaled variables $L$ and $s$, it means that we have the large $N$ behavior
\begin{equation}\\
E_{N,s}\left(L\right)\underset{N \to \infty}{\sim} \frac{399}{160} \sqrt{\frac{3}{\pi }} \Gamma \left(\frac{9}{4}\right) \, s \times N^{1/4}\quad   \hbox{for}\ 1\ll s\ll N^{1/4}\ .
\label{eq:EnsL}
\end{equation}
More generally, we find by expanding $\mathcal{F}(S,a,\lambda)$ at small $S$ the linear behavior:
\begin{equation}
\begin{split}
&E_{N,\left\lfloor S\, N^{1/4} \right\rfloor}\left(1-e^{-\lambda\, L/N^{1/2}}\right) \underset{N \to \infty}{\sim}
\psi(\lambda)\, S+O(S^2)\\
&\psi(\lambda)=\frac{133}{400}\sqrt{\frac{3}{2\pi}}\int_0^\infty\!\! \!\! dt\, t\, e^{-t^2}  \lambda ^{3/2} \left\{\!\left(15 t^2\!+\!4 \lambda ^2\right) \sqrt{\sqrt{1\!+\!\frac{4 t^2}{\lambda^2}}\!-\!1}\!-\!4 t\, \lambda\, \sqrt{\sqrt{1\!+\!\frac{4 t^2}{\lambda^2}}\!+\!1} \right\}\ . \\
\end{split}
\label{eq:psilambda}
\end{equation}
\begin{figure}
\begin{center}
\includegraphics[width=8cm]{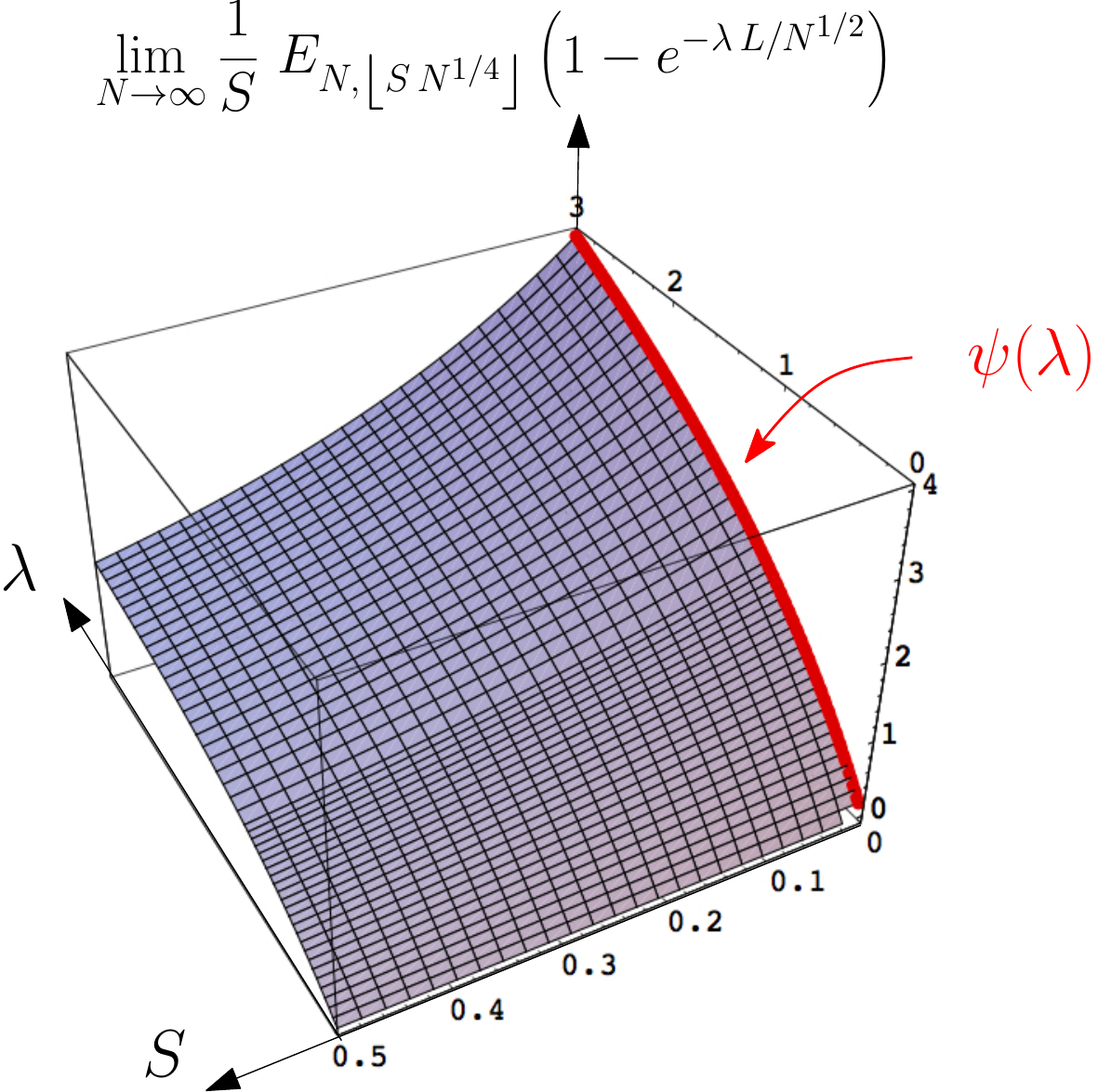}
\end{center}
\caption{A plot of the large $N$ limit of $\frac{1}{S}\ E_{N,\left\lfloor S\, N^{1/4} \right\rfloor}\left(1-e^{-\lambda\,  L/N^{1/2}}\right)$, which emphasizes
the small $S$ linear behavior \eqref{eq:psilambda}.}
\label{fig:checkpsiter}
\end{figure}
This behavior is emphasized in figure \ref{fig:checkpsiter} upon plotting 
$\lim\limits_{N\to\infty}\frac{1}{S}\ E_{N,\left\lfloor S\, N^{1/4} \right\rfloor}\left(1-e^{-\lambda\,  L/N^{1/2}}\right)$ as a function of $S$ and $\lambda$.

As for the large $S$ behavior of $\rho(S)$ and $\mu(S)$, a simple saddle point estimate of the corresponding contour integrals (in the $a$ variable) shows that
\begin{equation*}
\begin{split}
\rho(S) &\underset{S \to \infty}{\sim} 48\times 2^{1/6} 3^{5/6}\, S^{5/3}\,  e^{-3\, \left(\frac{3}{2}\right)^{2/3} S^{4/3}} \ ,\\
\mu(S) & \underset{S \to \infty}{\sim} 22\times \sqrt{6}\, S\,  e^{-3\, \left(\frac{3}{2}\right)^{2/3} S^{4/3}}\ ,  \\
\end{split}
\end{equation*}
so that we get the large $S$ behavior
\begin{equation}
  E_{N,\left\lfloor S\, N^{1/4} \right\rfloor}\left(L/N^{1/2}\right) \underset{N \to \infty}{\sim} \frac{11 }{12\times  2^{2/3}\, 3^{1/3}}\ \frac{1}{S^{2/3}}+O\left(\frac{1}{S^2}\right)\ .
\label{eq:largeS}
 \end{equation}
This limiting behavior is plotted in figure \ref{fig:aveLvsS}.
 
\section{Infinite maps: from scaling functions to the local limit} 

\subsection{Extracting the local limit from scaling functions}
The \emph{local limit} corresponds again to bi-pointed planar quadrangulations with a fixed volume $N$ and a fixed distance 
$2s$ between the two marked points, but in a regime where we now take the limit $N\to\infty$ while \emph{keeping $s$ finite}. 
It is now a general statement that we may describe properties of the local limit \emph{when $s$}, although remaining finite, \emph{becomes large}
in terms of the scaling functions that we already computed. Let us discuss this statement in more details.

The generating function $F_s(g,z)$ is singular when $g$ tends to its critical value $1/12$. For $g\to 1/12^{-}$, we have an expansion
of the form\footnote{It is easily checked by expanding \eqref{eq:eqforXst} in powers of $(1-12\, g)^{1/4}$ that $X_{s,t}(g,z)$ has an expansion 
which involves only even powers of $(1-12\, g)^{1/4}$, i.e.\ half-integer powers of $(1-12\, g)$, with moreover no term of order $(1-12\, g)^{1/2}$.
This property implies immediately a similar property for $F_s(g,z)=\Delta_s\Delta_t\log(X_{s,t}(g,z))\vert_{t=s}$. The
fact that $\mathfrak{f}_1(s,z)=0$ may also be understood by noting that $\mathfrak{f}_1(s,1)=0$, as easily checked from the explicit expression \eqref{eq:exactFgone},
and that $0\leq [g^N] F_s(g,z)\leq [g^N] F_s(g,1)$ for $0\leq z \leq 1$. Having a non-vanishing 
value of the leading singularity coefficient $\mathfrak{f}_1(s,z)$ when $0\leq z<1$ would indeed contradict this latter bound as it would
imply that $[g^N] F_s(g,z)$ grows faster than $[g^N] F_s(g,1)$ at large $N$.}
\begin{equation}
F_s(g,z)=\sum_{i\geq 0} \mathfrak{f}_i(s,z)\, (1-12\, g)^{i/2}
\label{eq:Fexpand}
\end{equation}
involving half-integer powers of $(1-12\, g)$, and with $\mathfrak{f}_1(s,z)=0$. 

Using the identification $(1-12\, g)^{1/2}=(a/\sqrt{6})^2 \epsilon^2$ whenever $g=G(a,\epsilon)$, 
we may then connect the scaling function $\mathcal{F}(S,a,\lambda)$ to the above expansion and write
\begin{equation*}
\begin{split}
\mathcal{F}(S,a,\lambda)&=\lim_{\epsilon\to 0}\,  \frac{1}{\epsilon^3} \, F_{\left\lfloor S/\epsilon\right\rfloor}\left(G(a,\epsilon),Z(\lambda,\epsilon)\right)\\
&=\lim_{\epsilon\to 0} \, \sum_{i\geq 0} \left(\frac{a}{\sqrt{6}}\right)^{2i}\, \epsilon^{2i-3}\, \mathfrak{f}_i\left(\left\lfloor S/\epsilon\right\rfloor,Z(\lambda,\epsilon)\right)\\
\end{split}
\end{equation*}
with $G(a,\epsilon)$ and $Z(\lambda,\epsilon)$ as in \eqref{eq:gzscal}. Now $Z(\lambda,\epsilon)$ depends only on the product $\lambda\, \epsilon^2$ so that $\mathfrak{f}_i\left(\left\lfloor S/\epsilon\right\rfloor,Z(\lambda,\epsilon)\right)$ depends on only \emph{two} variables $S/\epsilon$ and $\lambda\, \epsilon^2$, or equivalently on $S/\epsilon$ and $(S/\epsilon)^2 \times \lambda\, \epsilon^2 = \lambda\, S^2$. The existence of a non-trivial finite
limit when $\epsilon\to 0$ implies that 
\begin{equation}
\mathfrak{f}_i\left(\left\lfloor S/\epsilon\right\rfloor,Z(\lambda,\epsilon)\right) \underset{\epsilon \to 0}{\sim} \left(\frac{S}{\epsilon}\right)^{2i-3} \varphi_i(\lambda\, S^2)
\label{eq:fiphi}
\end{equation}
for $i\neq 1$. Plugging back this expression in $\mathcal{F}(S,a,\lambda)$, we then have the direct identification
\begin{equation*}
\hskip -1.2cm \mathcal{F}(S,a,\lambda)=
 \sum_{i\geq 0} \left(\frac{a}{\sqrt{6}}\right)^{2i}\, S^{2i-3}\, \varphi_i(\lambda\, S^2)
 \quad \Leftrightarrow \quad \mathcal{F}\left(S,a,\frac{\tau}{S^2}\right)=
 \sum_{i\geq 0} \left(\frac{a}{\sqrt{6}}\right)^{2i}\, S^{2i-3}\, \varphi_i(\tau)
\end{equation*}
with $\varphi_1=0$. 
We may therefore read the functions $\varphi_i(\tau)$ from the expansion in $S$ of $\mathcal{F}(S,a,\tau/S^2)$ via
\begin{equation}
\varphi_i(\tau)=  \left(\frac{a}{\sqrt{6}}\right)^{-2i} [S^{2i-3}] \mathcal{F}\left(S,a,\frac{\tau}{S^2}\right) =  [S^{2i-3}]  \mathcal{F}\left(S,\sqrt{6},\frac{\tau}{S^2}\right)\ ,
\label{eq:phival}
\end{equation}
where we eventually set $a=\sqrt{6}$ since $\varphi_i(\tau)$ does not depend on $a$.
\vskip .3cm
Returning to the local limit, we expect that, for $N\to \infty$ and for finite $s$, the Vorono\"\i\ cell perimeter $L$ will remain finite.
Moreover, as we shall see, $L$ scales as $s^2$ at large $s$ in the sense that 
the limiting expected value
\begin{equation*}
\Phi_s(\omega)\equiv \lim_{N\to \infty} E_{N,s}\left(e^{-\omega L/s^{2}}\right)
\end{equation*}
has a non-trivial limit when $s\to \infty$.
Let us indeed see how to estimate $\Phi_s(\omega)$ at large $s$ from the above correspondence. From \eqref{eq:ensdef}, we
have 
\begin{equation*}
\Phi_s(\omega)= \lim_{N\to \infty} \frac{[g^N]F_s\left(g,e^{-\omega/s^{2}}\right)}{[g^N]F_s(g,1)}
\end{equation*}  
with a leading singularity of $F_s(g,z)$ when $g\to 1/12^{-}$ corresponding to the $i=3$ term in \eqref{eq:Fexpand} since the $i=0$ and $i=2$ terms are regular and
the $i=1$ term is missing. This immediately yields the large $N$ estimate
\begin{equation*}
[g^N]F_s(g,z)\underset{N \to \infty}{\sim} \frac{3}{4}  \frac{12^{N}}{\sqrt{\pi} N^{5/2}} \mathfrak{f}_3(s,z)
\end{equation*}
from which we deduce the identity
\begin{equation*}
\Phi_s(\omega)= \frac{\mathfrak{f}_3\left(s,e^{-\omega/s^{2}}\right)}{\mathfrak{f}_3(s,1)}\ .
\end{equation*}  
Its large $s$ behavior may be obtained as follows: 
setting $\lambda=\omega/S^2$ and $\epsilon=S/s$ in \eqref{eq:fiphi}, with $\epsilon\to 0$ when $S$ is fixed and $s\to \infty$, we have, for $i=3$, 
\begin{equation*}
\mathfrak{f}_3\left(s,Z(\omega/S^2,S/s)\right) \underset{s \to \infty}{\sim} s^{3} \varphi_3(\omega)
=s^{3}  [S^{3}]  \mathcal{F}\left(S,\sqrt{6},\frac{\omega}{S^2}\right)
\end{equation*}
from \eqref{eq:phival}.
At large $s$, $Z(\omega/S^2,S/s)\sim e^{-\omega/s^{2}}$ and
we arrive immediately at the desired asymptotic expression
\begin{equation}
\Phi_s(\omega)\underset{s \to \infty}{\sim} \frac{[S^{3}]  \mathcal{F}\left(S,\sqrt{6},\displaystyle{\frac{\omega}{S^2}}\right)}{[S^3]\mathcal{F}(S,\sqrt{6},0)}
\label{eq:Phiomega}
\end{equation}
in terms of the scaling function $\mathcal{F}(S,a,\lambda)$ only.

\subsection{Explicit expressions and plots}
From the explicit expression \eqref{eq:Hval}--\eqref{eq:Fval}, it is a rather straightforward exercise to expand 
$ \mathcal{F}\left(S,\sqrt{6},\omega/S^2\right)$ in powers of $S$ and to extract its $S^3$ coefficient. 
The only non-trivial part of the exercise comes from the hypergeometric function which enters
the formula \eqref{eq:Hval} for $\widetilde{H}(\sigma,A)$. Its expansion is discussed in details in Appendix A and gives rise to
terms involving the exponential integral function $E_1(r)$ defined by
\begin{equation*}
E_1(r)\equiv \int_r^\infty \frac{dt}{t}\, e^{-t}\ .
\end{equation*}
The net result of the small $S$ expansion\footnote{As expected, the small $S$ expansion of $ \mathcal{F}\left(S,\sqrt{6},\omega/S^2\right)$
involves only odd powers of $S$, starts with a term of order $1/S^3$, has no term of order $1/S$ so that the next term is of order $S$, followed
by the desired term of order $S^3$.} is the following fully explicit expression:
\begin{equation*}
\begin{split}
& \hskip -1.2cm [S^{3}]  \mathcal{F}\left(S,\sqrt{6},\displaystyle{\frac{\omega}{S^2}}\right)=
\frac{1}{9800\, \left(\Pi(q)\right)^3} \left(\frac{\Lambda (q)}{\Pi (q)}+6\, (6q)^3\, \Omega(q)\, e^{6 q}\,  E_1(6 q) \right)\ ,\quad q=\sqrt{\frac{\omega}{3}}\ ,\\
& \Pi(q)=5+15 q+18 q^2+9 q^3 \ , \\
&\Omega(q)= 1925+5775 q+6930 q^2+3465 q^3+2700 q^6+2754 q^7+972 q^8 \ , \\
&\Lambda(q)= 14000000+110394375 q+394807500 q^2+841185000 q^3+1180872000 q^4 \\
&\qquad +1137241620 q^5+754502040 q^6+334611972 q^7+89532864 q^8
\\ & \qquad +3093876 q^9-17740944 q^{10}-17688456 q^{11}-8817984
   q^{12}-1889568 q^{13}\ .\\
\end{split}
\end{equation*}
For $\omega\to 0$, this expression tends to $16/7$ so that we obtain eventually the desired result
\begin{equation}
\Phi_s(\omega)\underset{s \to \infty}{\sim}
\frac{1}{22400\, \left(\Pi(q)\right)^3} \left(\frac{\Lambda (q)}{\Pi (q)}+6\, (6q)^3\, \Omega(q)\, e^{6 q}\,  E_1(6 q) \right)\ ,\quad q=\sqrt{\frac{\omega}{3}}\ .
\label{eq:valPhiomega}
\end{equation}
\begin{figure}
\begin{center}
\includegraphics[width=8cm]{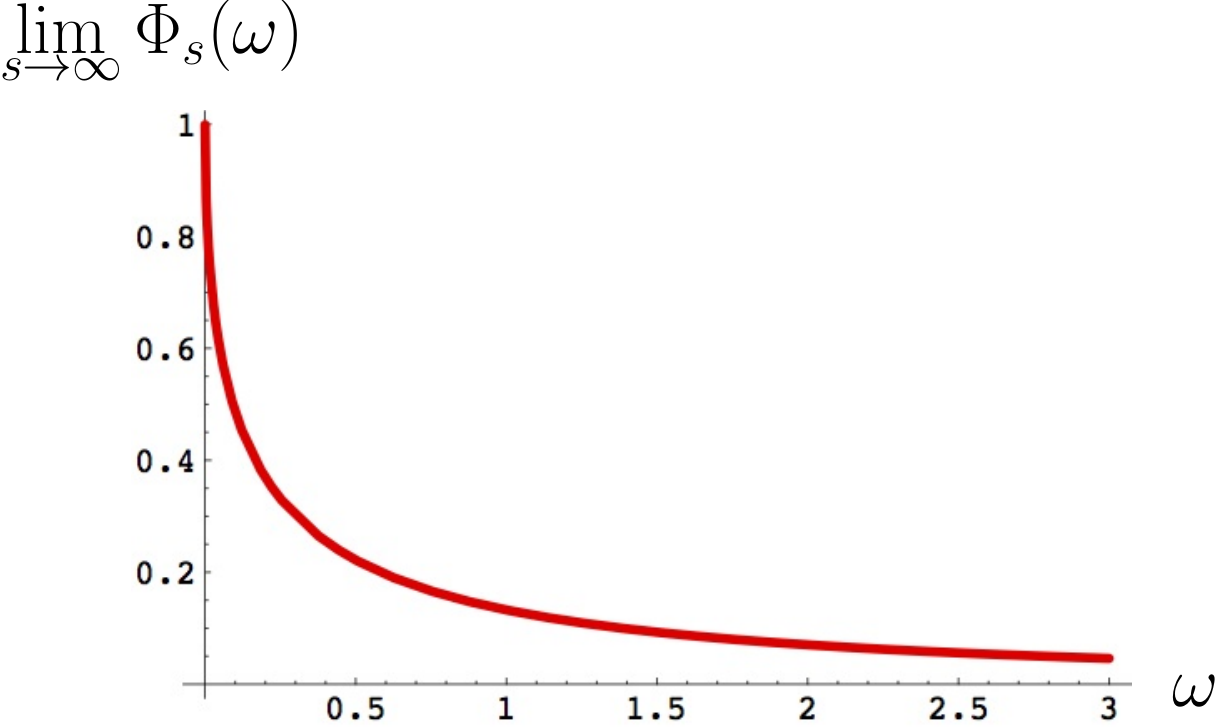}
\end{center}
\caption{A plot of the function $\Phi_s(\omega)$ at large $s$, as given by \eqref{eq:valPhiomega}. }
\label{fig:Phivsomega}
\end{figure}
This function is plotted in figure \ref{fig:Phivsomega} for illustration. 
At large $\omega$, it has the expansion
\begin{equation}
\lim_{s\to\infty} \Phi_s(\omega) = 
\frac{1}{4} \sqrt{3}\, \frac{1}{\omega^{3/2}}-\frac{15}{16} \sqrt{3}\, \frac{1}{\omega^{5/2}}+\frac{39}{16}\, \frac{1}{\omega^{3}}
+O\left(\frac{1}{\omega^{7/2}}\right)\ .
\label{eq:largeomega}
\end{equation}
At small $\omega$, its expansion reads
\begin{equation}
\lim_{s\to\infty} \Phi_s(\omega) = 1-\frac{4389 \sqrt{3}}{3200}\, \sqrt{\omega} +\frac{1713}{560}\,  \omega +O\left(\omega^{3/2}\log(\omega)\right)\ .
\label{eq:smallomega}
\end{equation}
In particular, due to the $\sqrt{\omega}$ singular term, the expected value of $L/s^2$ in the local limit, given by
\begin{equation*}
\lim_{N\to \infty} E_{N,s}\left(L/s^{2}\right)=-\frac{d}{d\omega} \Phi_s(\omega)\Big\vert_{\omega=0}
\end{equation*}
is \emph{infinite} for large $s$. This is consistent with the estimate \eqref{eq:EnsL} of the scaling limit which states that
\begin{equation*}
E_{N,s}\left(L/s^2\right)\underset{N \to \infty}{\sim} \frac{399}{160} \sqrt{\frac{3}{\pi }} \Gamma \left(\frac{9}{4}\right) \times \frac{N^{1/4}}{s} \quad   \hbox{for}\ 1\ll s\ll N^{1/4}
\end{equation*}
and therefore implies that $E_{N,s}\left(L/s^2\right)$ diverges as $N^{1/4}$ when $N\to \infty$ at fixed finite (although large) $s$.

The inverse Laplace transform (with respect to $\omega$) of $\Phi_s(\omega)$  is nothing but the \emph{probability density} 
$\mathcal{P}_s(\ell)$ for the
quantity 
\begin{equation*}
\ell\equiv \frac{L}{s^2}\ .
\end{equation*}
\begin{figure}
\begin{center}
\includegraphics[width=8cm]{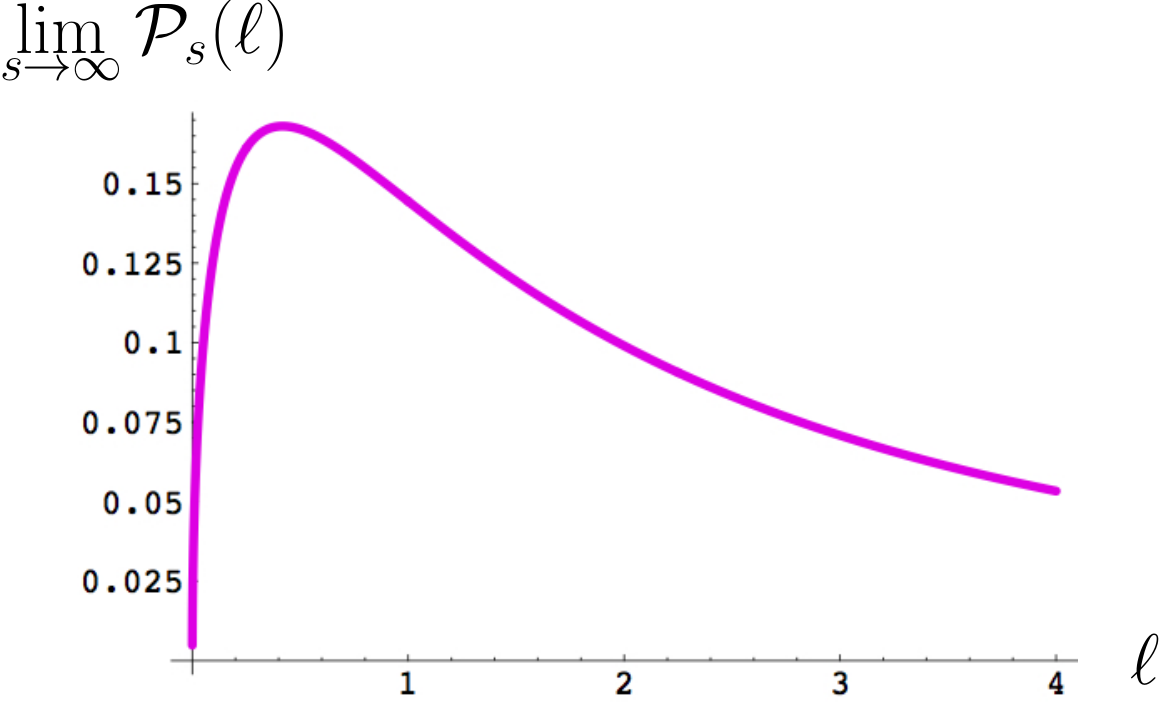}
\end{center}
\caption{A plot of the function probability density $\mathcal{P}_s(\ell)$ at large $s$, as obtained by a numerical inverse Laplace transform
of $\Phi_s(\omega)$. Its large $\ell$ and small $\ell$ behaviors are given by \eqref{eq:largeell} and \eqref{eq:smallell} respectively.}
\label{fig:ProbaL}
\end{figure}
It has a non-trivial limit at large $s$, which may be evaluated numerically from the explicit expression \eqref{eq:valPhiomega} via some appropriate
inverse Laplace transform numerical tool \cite{NILT,Ava,AVb}. The resulting plot is shown in figure \ref{fig:ProbaL}. For large $\ell$, we read from the small $\omega$ behavior
\eqref{eq:smallomega} that
\begin{equation}
\lim_{s\to\infty}\mathcal{P}_s(\ell)\underset{\ell \to \infty}{\sim} \frac{4389}{6400}\, \sqrt{\frac{3}{\pi}}\, \frac{1}{\ell^{3/2}}\ .
\label{eq:largeell}
\end{equation}
It particular, this behavior implies that all the positive moments of this law are infinite.
For small $\ell$, we read from the large $\omega$ behavior
\eqref{eq:largeomega} that
\begin{equation}
\lim_{s\to\infty}\mathcal{P}_s(\ell)\underset{\ell \to 0}{\sim} \frac{1}{2}\, \sqrt{\frac{3}{\pi}}\, \sqrt{\ell}\ .
\label{eq:smallell}
\end{equation}

\section{conclusion} 
In this paper, we presented a number of exact laws for the Vorono\"\i\ cell perimeter in large bi-pointed quadrangulations,
both in the scaling and in the local limit. Even though the explicit expressions of these laws, when available, are quite involved,
these formulas are expected to be universal and characteristic of the Brownian map (scaling limit) or the Brownian plane 
(local limit at large distance). In other words, the very same expressions would be obtained in the same regimes, up to some possible
 global rescalings 
(depending on the precise definition of distances, lengths and volumes), for all families of map in the universality
class of the so-called pure gravity, including in particular maps with faces of arbitrary but bounded degrees.

In the local limit $N\to \infty$ and $s$ finite, it was shown in \cite{G17b} that  exactly one of the two Vorono\"\i\ cells
becomes infinite while the other remains finite with a volume $V$ of order $s^4$. The Vorono\"\i\ cell perimeter 
$L$, which scales as $s^2$, is thus the length of the frontier between a finite and an infinite cell. It is interesting
to note that we found here 
that all the moments of the law for $L/s^2$ are infinite, in the same way as it was found in \cite{G17b} that all the moments of the law
for $V/s^4$ diverge.

Note finally that, if we do not fix the distance between the two marked vertices, we expect at large $N$ that the typical distance will be of order 
$N^{1/4}$ and we can therefore use the results of the scaling limit to state that the expected value
 of the Vorono\"\i\ cell perimeter in large bi-pointed quadrangulations having their two marked vertices at arbitrary even distance from each other
 will be asymptotically equal to $\gamma N^{1/2}$ with
 \begin{equation*}
 \gamma=\int_0^\infty
 dS\, \mu(S)=\frac{12\, \sqrt{\pi}}{35}\ .
 \end{equation*}
\section*{Acknowledgements} 
The author thanks Michel Bauer for enlightening discussions on various technical points in the computations. 
The author acknowledges the support of the grant ANR-14-CE25-0014 (ANR GRAAL).

\appendix
\section{Expansion of the hypergeometric function ${}\ _2F_1(1,2A+2,2A+3,\sigma)$} 
We are interested in computing the small $S$ expansion of $\mathcal{F}(S,\sqrt{6},\omega/S^2)$ up to order $S^3$.
To this end, we need from \eqref{eq:Hval} and \eqref{eq:Fval} the small $S$ expansion of the quantity
\begin{equation*}
\Xi(S,\omega)\equiv {}\ _2F_1(1,2A+2,2A+3,\sigma)\quad \hbox{with}\ \sigma=e^{-\sqrt{6}S}\ \hbox{and}\ A=\sqrt{1+\frac{\omega}{2S^2}}\ .
\end{equation*}
A first naive approach would consist in first expanding the term of order $\sigma^n$ in ${}\ _2F_1(1,2A+2,2A+3,\sigma)$, as 
read in \eqref{eq:def2F1}, for $A=\sqrt{1+\omega/(2S^2)}$, namely 
\begin{equation*}
\frac{2A+2}{n+2A+2}= 1-\frac{n}{\sqrt{2\omega }}\, S+\frac{n (2+n)}{2 \omega }\, S^2+\cdots
\end{equation*}
and in summing over $n$ to get the small $S$ expansion
\begin{equation*}
\hskip -1.2cm {}\ _2F_1(1,2A+2,2A+3,\sigma)=\frac{1}{1-\sigma}-\frac{1}{\sqrt{2\omega }}\, \frac{\sigma}{(1-\sigma)^2}\, S+\frac{1}{2 \omega }\, 
\frac{3\sigma-\sigma^2}{(1-\sigma)^3}\, S^2+\cdots\ .
\end{equation*}
Setting then, in a second step, $\sigma=e^{-\sqrt{6}S}$ and expanding at small $S$, we deduce
\begin{equation}
\Xi(S,\omega)=\left(\frac{1}{\sqrt{6} S}+\cdots\right)-\frac{1}{\sqrt{2\omega }}\, \left(\frac{1}{6 S}+\cdots\right)+\frac{1}{2 \omega }\, 
\left(\frac{1}{3\sqrt{6}S}+\cdots\right)+\cdots
\label{eq:naive}
\end{equation}
where it is apparent that \emph{all the terms of the expansion contribute} in fact \emph{ to the same leading order} $1/S$. This naive approach therefore 
fails and a proper resummation is required.  
To perform this resummation and properly get the desired small $S$ expansion, we may instead rely on the differential equation satisfied by 
${}\ _2F_1(1,2A+2,2A+3,\sigma)$, namely
\begin{equation*}
\hskip -1.2cm \sigma\, \frac{d}{d\sigma}{}\ _2F_1(1,2A+2,2A+3,\sigma) +(2A+2){}\ _2F_1(1,2A+2,2A+3,\sigma)-\frac{2A+2}{1-\sigma}=0\ ,
\end{equation*}
which immediately translates into the following equation for $\Xi(S,\omega)$:
\begin{equation}
\hskip -.1cm \frac{\partial}{\partial S}\Xi(S,\omega)+2\frac{\omega}{S}\frac{\partial}{\partial\omega}\Xi(S,\omega)-2\sqrt{6}\, \left(\sqrt{1+\frac{\omega}{2S^2}}+1\right)\, \Xi(S,\omega)
+2\sqrt{6}\, \frac{\sqrt{1+\displaystyle{\frac{\omega}{2S^2}}}+1}{1-e^{-\sqrt{6}S}}=0\ .
\label{eq:diffeq}
\end{equation}
From the previous discussion, $\Xi(S,\omega)$ is expected to have a small $S$ leading term of 
order $1/S$, in which case all the different terms in the differential equation above are precisely of the same order $1/S^2$.
Writing therefore
\begin{equation*}
\Xi(S,\omega)=\xi_{-1}\left(2\sqrt{3\omega}\right)\, \frac{1}{S/\sqrt{6}}+O(1)\ ,
\end{equation*}
where the multiplicative factors were introduced for convenience, and expanding \eqref{eq:diffeq} to leading order $1/S^2$ at small $S$
yields the following differential equation for $\xi_{-1}(r)$:
\begin{equation*}
6 r \, \frac{d}{dr}\xi_{-1}(r)-6(1+r)\, \xi_{-1}(r)+r=0\ ,
\end{equation*}
easily solved into
\begin{equation*}
\xi_{-1}(r)=k\, r\, e^r+\frac{1}{6}\, r\, e^r\, E_1(r)\ ,
\end{equation*}
where $k$ is some arbitrary constant and $E_1(r)$ the exponential integral function
\begin{equation*}
E_1(r)\equiv \int_r^\infty \frac{dt}{t}\, e^{-t}\ .
\end{equation*}
When $\omega\to \infty$, we have clearly $\Xi(S,\omega)\to 1/(1-e^{-\sqrt{6}\, S})=1/(\sqrt{6}S)+O(1)$ so that $\xi_{-1}(r)$ must have a \emph{finite large $r$ limit} equal to $1/6$. The existence of this finite limit at large $r$ fixes the integration constant $k$ to $k=0$, so that
\begin{equation*}
\xi_{-1}(r)=\frac{1}{6}\, r\, e^r\, E_1(r)\ .
\end{equation*}
Note that, from the ``naive" analysis, we expect moreover the following large $r$ expansion,
read from \eqref{eq:naive}:
\begin{equation*}
\xi_{-1}(r)=\frac{1}{6}\left(1-\frac{1}{r}+\frac{2}{r^2}+O\left(\frac{1}{r^3}\right)\right)\ .
\end{equation*}
These expansion is indeed recovered from the general property
\begin{equation*}
E_1(r)=\frac{e^{-r}}{r}\sum_{p=0}^P\frac{p!}{(-r)^p}+O\left(\frac{e^{-r}}{r^P}\right)\ .
\end{equation*}
To get the desired expansion of $\mathcal{F}(S,\sqrt{6},\omega/S^2)$ up to order $S^3$, it is easily checked from \eqref{eq:Hval} and \eqref{eq:Fval} 
that the small $S$ expansion of $\Xi(S,\omega)$ must be push up to order $S^5$. We therefore write
\begin{equation*}
\Xi(S,\omega)=\sum_{i=-1}^5 \xi_{i}\left(2\sqrt{3\omega}\right) (S/\sqrt{6})^i+O(S^6)\ .
\end{equation*}
Expanding \eqref{eq:diffeq} to increasing orders in $S$ gives rise to a sequence of differential equations for the successive 
functions $\xi_i(r)$, $i=0,1,\dots,5$ which are easily solved, one after the other, leading to
 \begin{equation*}
 \begin{split}
 \xi_0(r)&= \frac{1}{2} \Big(4\,  (1+r)\, e^r\, E_1(r)-3\Big)\ ,\\
 \xi_1(r)&= \frac{1}{2\, r}\Big(24\,  \left(1+3 r+r^2\right) \,  e^r\, E_1(r)-(47+23 r)\Big)\ ,\\
 \xi_2(r)&= \frac{6}{r^2} \Big(8\, r \left(6+6 r+r^2\right)\, e^r\, E_1(r)-(17+40 r+8 r^2)\Big)\ ,\\
 \xi_3(r)&= \frac{3}{10\, r^3} \Big(480 \left(-3+3 r+21 r^2+10 r^3+r^4\right)\, e^r\, E_1(r)\\ & \qquad \qquad \qquad \qquad 
 -(-2034+6246 r+4323 r^2+481 r^3)\Big)\ ,\\
 \xi_4(r)&= \frac{54}{5\, r^4} \Big(32\, r \left(-30+30 r+55 r^2+15 r^3+r^4\right)\,  e^r \, E_1(r)\\
 &\ \  \qquad \qquad \qquad-(-518-36 r+1343 r^2+448 r^3+32 r^4)\Big)\ ,\\
  \xi_5(r)&= \frac{9}{35\, r^5}\Big(2688 (45-45 r-135 r^2+150 r^3+120 r^4+21 r^5+r^6)\, e^r\, 
   E_1(r)\\&\ \ \qquad  -(154608-393240 r+177600 r^2+271468 r^3+53755 r^4+2687 r^5)\Big)\ .\\
 \end{split}
 \end{equation*}
This provides the desired small $S$ expansion of $\Xi(S,\omega)$ which, inserted in \eqref{eq:Hval}--\eqref{eq:Fval}, yields
the wanted expansion of $\mathcal{F}(S,\sqrt{6},\omega/S^2)$ up to order $S^3$.
\bibliographystyle{plain}
\bibliography{voronoilength}
\end{document}